\documentstyle[12pt]{article}
\newlength{\abstractwidth}
\renewcommand{\thefootnote}{\fnsymbol{footnote}}
\renewcommand{\thanks}[1]{\footnote{#1}} 
\newcommand{\starttext}{ \setcounter{footnote}{0}
\renewcommand{\thefootnote}{\arabic{footnote}}}

\newcommand{\be}{\begin{equation}}
\newcommand{\bea}{\begin{eqnarray}}
\newcommand{\eea}{\end{eqnarray}} \newcommand{\ee}{\end{equation}}
 
 \def\ba{\begin{eqnarray}}
\def\ea{\end{eqnarray}}

\def\o{\omega}

\def\tr{{\rm tr}}
\def\det{{\rm det}}

\def\log{\,{\rm log}\,}
\def\exp{\,{\rm exp}\,}

\def\o{\omega}

\def\o{\omega}

\def\ge{\geq}
\def\le{\leq}

\def\p{\partial}

\def\[{{\bf [}}
\def\]{{\bf ]}}

\def\ddbar{i\p\bar\p}

\def\mathbb{\bf}

\def\eqref{\ref}

\newcommand{\neweqref}[1]{(\ref{#1})}

\begin{document}
\starttext \baselineskip=18pt \setcounter{footnote}{0}
\newtheorem{theorem}{Theorem}
\newtheorem{lemma}{Lemma}
\newtheorem{corollary}{Corollary}
\newtheorem{definition}{Definition}
\newtheorem{conjecture}{Conjecture}
\newtheorem{proposition}{Proposition}


\begin{center}
{\Large \bf ON $L^\infty$ ESTIMATES FOR FULLY NONLINEAR PARTIAL DIFFERENTIAL EQUATIONS
\footnote{Work supported in part by the National Science Foundation under grants DMS-1855947 and DMS-22-03273.}}

\medskip
\centerline{Bin Guo and  Duong H. Phong}

\medskip

\begin{abstract}

{\footnotesize Sharp $L^\infty$ estimates are obtained for general classes of fully non-linear PDE's on non-K\"ahler manifolds, complementing the theory developed earlier by the authors in joint work with F. Tong for the K\"ahler case. The key idea is still a comparison with an auxiliary Monge-Amp\`ere equation, but this time on a ball with Dirichlet boundary conditions, so that it always admits a unique solution. The method applies not just to compact Hermitian manifolds, but also to the Dirichlet problem, to open manifolds with a positive lower bound on their injectivity radii, to $(n-1)$ form equations, and even to non-integrable almost-complex or symplectic manifolds. It is the first method applicable in any generality to large classes of non-linear equations, and it usually improves on other methods when they happen to be available for specific equations. }

\end{abstract}

\end{center}

\baselineskip=15pt
\setcounter{equation}{0}
\setcounter{footnote}{0}

\section{Introduction}
\setcounter{equation}{0}

A priori estimates are essential to partial differential equations, and of these, $L^\infty$ estimates are of particular importance, as they are usually also required for other estimates. In particular, the modern theory of partial differential equations can be argued to have started with the De Giorgi-Nash-Moser theory of $L^\infty$ estimates for linear equations in divergence form. The importance of $L^\infty$ estimates for non-linear equations was recognized early on in the case of complex Monge-Amp\`ere equations on K\"ahler manifolds, where they were first obtained by Yau \cite{Y} in his celebrated solution of the Calabi conjecture.
A sharp version was subsequently obtained by Kolodziej \cite{K} using pluripotential theory. Very recently, an alternative proof of Kolodziej's estimates using only PDE methods was provided by the authors in joint work with F. Tong \cite{GPT}. This proof also gave new and sharp estimates for general classes of non-linear equations on K\"ahler manifolds, thus providing a first step towards a De Giorgi-Nash-Moser theory for fully non-linear equations in complex geometry.

\medskip
The goal of the present paper is to obtain sharp $L^\infty$ estimates for fully non-linear equations in non-K\"ahler settings. There is a strong incentive for this, as non-K\"ahler settings appear increasingly frequently in complex geometry, symplectic geometry, and theoretical physics (see e.g. \cite{FR, FY, GF, P, Po} for surveys). Extensions from K\"ahler to non-K\"ahler settings are typically technically much more complicated, as new torsion terms can proliferate and integration by parts is unwieldy. However, there are usually considerable new conceptual difficulties as well, because of a lack of cohomological constraints. In particular, because of this, a basic idea from \cite{GPT}, namely a comparison with an auxiliary global complex Monge-Amp\`ere equation, runs immediately into difficulties: in the absence of suitable cohomological constraints, the auxiliary equation may not be solvable.

\medskip
In this paper, we shall show that comparisons with auxiliary Monge-Amp\`ere equations are still of powerful assistance, but the auxiliary equations should not be defined globally, but only on an open ball with Dirichlet conditions. The existence and regularity of the auxiliary problem which we then need is a classic result of Caffarelli, Kohn, Nirenberg, and Spruck \cite{CKNS}, which plays the same role for the Dirichlet problem as Yau's existence theorem \cite{Y} did in \cite{GPT} for the compact K\"ahler case. An important point is that the desired global estimates can still be obtained from the local estimates obtained in this manner near a minimum. This requires an adaptation of the original arguments of Blocki \cite{B11}, combined with weak Harnack inequalities for supersolutions of linear elliptic equations \cite{GT}, as well as an estimate of Kolodziej \cite{K} which generalizes the one-dimensional estimate of Brezis-Merle \cite{BM}. 

\medskip
The resulting method for $L^\infty$ estimates, using only comparisons with a local equation, turns out to be both flexible and powerful. It bypasses the complicated torsion terms and integration by parts plaguing non-K\"ahler geometry mentioned above, and it applies to whole new classes of fully non-linear equations, not just on compact Hermitian manifolds, but also to the Dirichlet problem, to open manifolds, to $(n-1)$-form equations, and to equations on almost-complex manifolds. The extension to almost-complex manifolds, using auxiliary real Monge-Amp\`ere equations, is in particular of considerable interest, 
as non-integrable almost-complex manifolds have revealed themselves to be effective intermediaries between complex geometry and symplectic geometry (see e.g. \cite{Hit, D, TWY, FPPZ}). On the other hand, if we assume the K\"ahler condition, then comparisons with a global equation as in \cite{GPT} are possible, and we obtain $L^\infty$ estimates which have the great advantage of remaining uniform over classes of K\"ahler metrics possibly degenerating to the boundary of the K\"ahler cone. Such a uniformity is essential for many applications to algebraic geometry (see e.g. \cite{EGZ, DP, GPSS}). Thus the present paper and \cite{GPT} address different geometric requirements in complex geometry, and together, they are a natural candidate for what may be viewed as a De Giorgi-Nash-Moser theory for non-linear equations in this area.

\medskip
We describe now more fully our results. 
Let $(X,\omega)$ be a connected complex manifold of complex dimension $n$,  and $\omega=\sum_{jk}ig_{\bar kj} dz^j\wedge d\bar z^k$ be a Hermitian metric. Let $f: \Gamma \subset {\mathbb R}^n\to {\mathbf R}_+$ be 
a function satisfying the following conditions: 

 (1) $\Gamma\subset {\mathbb R}^n$ is a symmetric cone with 
\begin{equation}\label{eqn:cone}
\Gamma_n\subset \Gamma \subset \Gamma_1
\end{equation}
where $\Gamma_k$ is the cone of vectors $\lambda$ with $\sigma_j(\lambda)>0$ for $1\leq j\leq k$. Here $\sigma_k(\lambda)$ is the $j$-th symmetric polynomial in $\lambda$. In particular, $\Gamma_1$ is the half-space defined by $\lambda_1+\cdots+\lambda_n>0$, and $\Gamma_n$ is the first octant, defined by $\lambda_j>0$ for $1\leq j\leq n$.

(2) $f(\lambda)$ is symmetric in $\lambda = (\lambda_1,\ldots, \lambda_n)\in \Gamma$ and it is homogeneous of degree one;

 (3) $\frac{\partial f}{\partial \lambda_j}>0$ for each $j=1,\ldots, n$ and $\lambda\in \Gamma$;

 (4) There is a constant $\gamma>0$ such that 
\begin{equation}\label{eqn:structure}
\prod_{j=1}^n \frac{\partial f}{\partial \lambda_j}\ge \gamma,\quad \forall \lambda\in \Gamma.
\end{equation}

\medskip

Familiar examples of such nonlinear operators $f(\lambda)$ include
the complex Monge-Amp\`ere operator $f(\lambda[h_\varphi]) = 
(\prod_{j=1}^n\lambda_j)^{1\over n}$, $\Gamma = \Gamma_n$;
the complex Hessian operator $f(\lambda[h_\varphi]) = (\sigma_k(\lambda))^{1\over k}$,
$\Gamma = \Gamma_k$ where $\sigma_k(\lambda)$ is the $k$-th symmetric polynomial in the eigenvalues $\lambda_j$ \cite{CNS,DK2,TW};
and the $p$-Monge-Amp\`ere operator of Harvey and Lawson \cite{HL, HLnew},
$f(\lambda) = \Big( \prod_{I} \lambda_I\Big)^{\frac{n!}{(n-p)!p!}}$,
 where $I$ runs over all distinct multi-indices $1\le {i_1}<\cdots < {i_p}\le n$, $\lambda_I = \lambda_{i_1} + \cdots + \lambda_{i_p}$, and $\Gamma$ is the cone defined by $\lambda_I>0$ for all  $p$-indices $I$.
Clearly, the above conditions (1-4) are preserved on the intersection of the individual cones by finite linear combinations with positive coefficients of operators $f(\lambda)$ satisfying them. Remarkably, it has recently been shown by Harvey and Lawson \cite{HL22} that the set of operators $f(\lambda)$ satisfying the conditions (1-4) is quite large, and actually includes all invariant Garding-Dirichlet operators.
 
 \medskip
 
For any $C^2$ function $\varphi: X\to {\mathbb R}$, we set $\omega_\varphi = \omega + \ddbar \varphi$ and let $h_\varphi: T^{1,0}X\to T^{1,0}X$ be the endomorphism defined by $\omega_\varphi$ relative to $\omega$, i.e. in local coordinates, $(h_\varphi)^j{}_i= g^{j\bar k} (g_\varphi)_{\bar k i}$, where $g_{\bar j i}$ denotes the components of $\omega$ and $(g^{i\bar j})$ its inverse. We denote by $\lambda[h_\varphi]$ the (unordered) vector of eigenvalues of $h_\varphi$. 
Given a smooth function $F$ on $X$, we consider the fully nonlinear partial differential equation
\begin{equation}\label{eqn:main}
f(\lambda[h_\varphi]) = e^{F},
\mbox{ and } \lambda[h_\varphi]\in \Gamma. 
\end{equation}
Such equations were introduced systematically in the foundational paper by Caffarelli, Nirenberg, and Spruck \cite{CNS}, and occur now regularly in geometric analysis, e.g. \cite{GG, TW}.

\medskip

\begin{theorem}\label{thm:main1}\label{thm:main}
Assume that $(X,\o)$ is a compact Hermitian manifold without boundary, and consider the equation \neweqref{eqn:main} with the operator $f(\lambda)$ satisfying the conditions (1-4).
Then for any $p>n$ and any $C^2$ solution $\varphi$ to \neweqref{eqn:main}, normalized to satisfy
$\sup_X \varphi = 0$, we have 
\begin{equation}\label{eqn:C0}
{\sup}_X|\varphi|\le C,
\end{equation}
for some constant $C>0$ depending only on $X,\o, n, p, \gamma$, and $\| e^{nF}\|_{L^1(\log L)^p}$. Here the norm $\| e^{nF}\|_{L^1(\log L)^p}$ is defined by 
$$\| e^{nF}\|_{L^1(\log L)^p} = \int_X (1+ |F|^p) e^{nF} \omega^n$$ and is sometimes called the $p$-th Nash entropy of the function $e^{nF}$. 

\end{theorem}

\medskip
We survey briefly the previous known results related to this theorem. The case $(X,\o)$ K\"ahler was established in \cite{GPT}, including even a stronger version allowing degenerations of the underlying K\"ahler class $\o$ generalizing the case of the Monge-Amp\`ere equation obtained in \cite{EGZ, DP}. So we concentrate on the case of $(X,\o)$ non-K\"ahler, which is the focus of the present paper. In this case, for the Monge-Amp\`ere equation, $L^\infty$ estimates in the non-K\"ahler case were first obtained by Tosatti-Weinkove \cite{ToWe}, building on earlier works of Cherrier \cite{Ch}, and a pointwise upper bound for $e^{nF}$ was needed.
The stronger version requiring only weaker entropy bounds as described in Theorem \ref{thm:main1} was first obtained by Dinew and Kolodziej \cite{Di,DK1,DK2,K}, using a non-trivial extension of pluripotential theory to the Hermitian setting. More recently a new approach using the theory of envelopes has been developed by Guedj and Lu \cite{GL}. The proof of Theorem \ref{thm:main1} is arguably simpler than in all these approaches. More importantly, it is a PDE-based proof, and can apply to many more equations. The estimates in Theorem \ref{thm:main1} appear to be the first $L^\infty$ estimates obtained in any generality for fully non-linear equations on compact Hermitian manifolds.

\medskip
It may be worth noting an intriguing consequence of Theorem \ref{thm:main1}, or rather of its proof, as we shall see later in Section \S 3.
Given a $C^2$ function $\varphi$ such that $\omega_\varphi = \omega + \ddbar \varphi>0$ defines a Hermitian metric on $X$, in general the following identity does not hold
$$\int_X (\omega + \ddbar \varphi)^n = \int_X \omega^n,$$
unless extra conditions  \cite{Di} are put on $\omega$. So it is interesting to estimate the lower/upper bound of the volume of $(X,\omega_\varphi),$ i.e. the integral $\int_X \omega_\varphi^n$.
To this end, we define the relative volume $e^{nF}$ of $\omega_\varphi$ with respect to $\omega$ by $e^{nF} = \omega_\varphi^n/\omega^n$. Then $\varphi$ can be viewed as the solution of a complex Monge-Amp\`ere equation:
\begin{equation}\label{eqn:MA}
(\omega+ \ddbar \varphi)^n = e^{nF}\omega^n.
\end{equation}
Normalizing $\varphi$ so that $\sup_X \varphi = 0$, we can apply Theorem \ref{thm:main1} and its proof and obtain:

\begin{corollary}\label{cor:main}
Given $p>n$ and $K>0$, suppose $e^{nF}$ satisfies $\int_X e^{nF} |F|^p \omega^n \le K$, then there is a constant $c>0$ which depends on $n, \omega, p$ and $K$ such that 
\begin{equation}\label{eqn:vol lower}\int_X (\omega + \ddbar \varphi)^n = \int_X e^{nF} \omega^n \ge c.\end{equation}
\end{corollary}

Thus an upper bound on the $p$-th entropy of relative volume $e^{nF}$ implies a positive lower bound of the $L^1$-norm of $e^{nF}$. This is trivial if $\omega$ is K\"ahler, but not in general. We remark that  if $X$ admits a closed real $(1,1)$-form $\beta$ whose Bott-Chern class $[\beta]\in H^{1,1}_{BC}(X, {\mathbb C})$ is nef and big (i.e. $X$ is of {\em Fujiki class ${\mathcal C}$}), by a weak transcendental Morse inequality of Demailly  \cite{D}, equation \neweqref{eqn:vol lower} holds even without the assumption on $\int_X e^{nF} |F|^p \omega^n $. This is based on an observation of Tosatti \cite{T}. For the convenience of readers, we provide a proof (see Lemma \ref{lemma Demailly} in Section \S\ref{section 3 new}) of this fact using the arguments in \cite{D}.

\bigskip

The strategy of the proof of Theorem \ref{thm:main1} is to employ a {\em local} auxiliary complex Monge-Amp\`ere equations as in \cite{GPT} and a localized argument similar to that in \cite{B11}. As stressed above, this method is effective in a vast range of different geometric situations. We discuss some of these next.

\bigskip
{\bf The Dirichlet problem on Hermitian manifolds}

We consider first the case of the Dirichlet problem. Thus let $X$ be an open subset with smooth boundary in a larger Hermitian manifold $(\tilde X,\o)$.
Let $f(\lambda)$ be an operator satisfying the conditions (1-4), $\rho\in C^0(\p X)$, and consider the following Dirichlet boundary value problem on $X$:
\begin{equation}\label{eqn:DMA}
f(\lambda[h_\varphi]) = e^{F} \quad \mbox{in } X,\quad \lambda[h_\varphi]\in\Gamma,
\quad
\varphi = \rho \ {\rm on}\ \partial X.
\end{equation}

\begin{theorem}\label{thm:main2}
Given $p>n$, suppose $\varphi\in C^2(X)\cap C^0(\overline{ X})$ is a solution to the above Dirichlet problem \neweqref{eqn:DMA}. Then there is a constant $C>0$ which depends only on $X,\o,n,  p,\gamma$ and, $\| e^{nF}\|_{L^1(\log L)^p}$, such that 
\begin{equation}\label{eqn:D estimate}
{\rm sup}_X|\varphi| \le \| \rho\|_{L^\infty(\partial X)} + C(1+ \| \varphi\|_{L^1}^2).
\end{equation}
\end{theorem}

\medskip
We observe that, for complex Monge-Amp\`ere and Hessian equations, and if $\omega$ is K\"ahler, then $\| \varphi\|_{L^1}$ can be bounded in terms of the other parameters (see \S 4.6). More importantly, for general equations $f(\lambda)$, there are very few $L^\infty$ estimates, if any, which require only weak integral bounds on the right hand side $e^F$ such as in Theorem \ref{thm:main2}. For example, the $L^\infty$ estimates obtained by S.Y. Li in \cite{Li} are based on the construction of subsolutions, and as is usually the case with subsolutions, they require $L^\infty$ bounds on the right hand side.

\bigskip

{\bf Open Hermitian manifolds}

Next, we consider the case of a noncompact complex manifold $(X,\omega)$, where $\omega$ is a complete Hermitian metric on $X$. We shall make the following assumption:

\medskip
 
(A) There exists a constant $r_0>0$, so that, at each point $x\in X$, there exists a holomorphic coordinates chart $(U_x, z)$ centered at $x$ satisfying 
\begin{equation}\label{eqn:assumptions}\{|z|< 2 r_0\}\subset U_x, \quad \frac 12 \omega  \le \ddbar |z|^2 \le 2 \omega \mbox{ in }U_x.\end{equation}

\medskip

With the assumption \neweqref{eqn:assumptions}, we introduce two quantities associated to $\varphi$ and $e^F$, which we assume are bounded. Given a fixed $p>n$, we set 
\begin{equation}\label{eqn:varphi F}
K_\varphi = \sup_{x\in X} \int_{U_x} |\varphi| \omega^n, \quad K_F = \sup_{x\in X} \int_{U_x} e^{nF} |F|^p\omega^n.
\end{equation}

\begin{theorem}\label{thm:main3}
Let $(X,\o)$ be a noncompact Hermitian manifold satisfying the assumption (A). 
Let $\varphi\in L^\infty(X)\cap C^2(X)$ be a solution to \neweqref{eqn:DMA}.
Then there is a constant $C>0$ which depends on $n, p, \omega, \gamma$ and $K_\varphi, K_F$ such that 
$${\sup}_X|\varphi|\le C + {\rm liminf}_{z\to\infty}(-\varphi)_+$$
where $C$ is a constant depending on $X,\o,n, p, \gamma, r_0$, and $K_\varphi, K_F$.
\end{theorem}

\medskip

We remark that the assumption \neweqref{eqn:assumptions} is satisfied if the metric $\omega$ is K\"ahler, the Riemannian curvature is bounded and the injectivity radius of $(X,\omega)$ at each point is bounded from below by some positive number (see \cite{TY}). In particular, it holds for {\em asymptotically locally Euclidean} (ALE) K\"ahler manifolds, or under some asymptotic conditions on the (K\"ahler) metric $\omega$ near infinity, such as the {\em asymptotically conical} (AC) and {\em asymptotically cylindrical} (ACyl) conditions, which have been extensively studied in recent years. Typically the conditions which have been imposed in the literature on the open complete K\"ahler case are stronger than in the above theorem, see e.g. \cite{TY}. But more importantly, we are not aware of previous results for general fully non-linear equations on open complete Hermitian manifolds.

\bigskip
{\bf General $(n-1)$-form equations with gradient terms}

We show now how the local comparison method developed in this paper can apply to other second-order equations and not just  functions of the eigenvalues of the standard Hessian $h_\varphi$. Well-known examples of such equations arose in \cite{FWW, Po} and in the solution of the Gauduchon conjecture \cite{STW, TW, Po}, and they have been attracting considerable attention ever since \cite{GG, GGQ}. 

\smallskip
More specifically, 
let $(X,\omega)$ be a compact Hermitian manifold with $\partial X = \emptyset$ and $\omega_h$ be another Hermitian metric on $X$. Fix a smooth $(1,0)$-form $\theta$ on $X$, set
\bea
\chi[\varphi]=
i(\varphi_j\bar\theta_k+\theta_j\varphi_{\bar k})dz^j\wedge d\bar z^k,
\eea
and consider the endomorphism $\tilde h_\varphi$ of the tangent bundle $T^{1,0}(X)$ given by
\bea
\tilde h_\varphi
=
\omega^{-1}\cdot \big\{\omega_h+{1\over n-1}((\Delta_\o\varphi)\omega
-i\p\bar\p\varphi)+\chi[\varphi])\big\}
\eea
where 
$$\Delta_\omega \varphi = n\frac{\ddbar \varphi \wedge \omega^{n-1}}{\omega^n}$$ is the complex Laplacian of $\varphi$ with respect to $\omega$. We consider the fully non-linear equation

\begin{equation}\label{eqn:HLMA}
f(\lambda[\tilde h_\varphi]) = e^F,\quad \lambda[\tilde h_\varphi]\in \Gamma,
\end{equation} 
with $\varphi$ normalized by $\sup_X \varphi = 0$. We have then:

\begin{theorem}\label{thm:main4}
Given $p>2n$, there exists a constant $C>0$ which depends on $n,\gamma,  \omega, p, \theta, \omega_h$ and $\| e^{2 nF}\|_{L^1(\log L)^p(X,\omega)}$ such that 
$${\sup}_X |\varphi| \le C.$$
\end{theorem} 

\medskip
Previously, this equation had received most attention when $f(\lambda)=(\prod_{j=1}^n\lambda_j)^{1\over n}$, $\chi[\varphi]=0$, and was known in this case as the $(n-1)$-form Monge-Amp\`ere equation.
It was shown by Tosatti-Weinkove \cite{TW1} to be solvable up to a modification of $F$ by a suitable additive constant. This required a priori estimates of all orders, of which, as stressed in \cite{TW1}, the $L^\infty$ estimate is a major one.

In an earlier version of the present paper (arXiv:2204.12549), we had given a new and simpler proof of the $L^\infty$ estimate of Tosatti and Weinkove using our present method. Our proof was subsequently extended by Klemyatin, Liang, and Wang \cite{KLW} to more general operators $f(\lambda[\tilde h_\varphi])$ satisfying the structural conditions (1-4). 

The above Theorem \ref{thm:main4} shows that $L^\infty$ bounds can still be established even with the incorporation of gradient terms $\chi[\varphi]$ in the equation $f(\lambda[\tilde h_\varphi])$. Such terms are actually required for the solution of the Gauduchon conjecture \cite{STW}, and in general, they can lead to considerable additional difficulties in a priori estimates (see e.g. \cite{STW, GN}). Here, we show that these terms can readily be handled by our method, just by using as auxiliary equation the Dirichlet problem for a {\it real} Monge-Amp\`ere equation instead of the usual complex one. The only cost of using a real auxiliary equation is the constraint $p>2n$ in the $p$-Nash entropy, instead of the milder condition $p>n$ when $\chi=0$ (arXiv:2204.12549, \cite{KLW}). It is not known whether $p>n$ would suffice when $\chi\not=0$.

\bigskip
{\bf Fully non-linear equations on almost-complex manifolds}

While the importance of almost-complex structures has been known for a long-time, the structures which are not integrable had been comparatively neglected. This situation changed drastically in the last twenty years or so, as basic geometric structures, notably forms in lower dimensions or symplectic structures, have led to non-integrable almost-complex structures. We refer to \cite{Hit, D, TWY} for more detailed discussions. Here we would like to point out that the methods of this paper readily adapt to produce $L^\infty$ estimates for general classes of non-linear equations, even on almost-complex manifolds.

\medskip
Let $X$ be a compact manifold without boundary of real dimension $2n$, equipped with an almost complex structure $J$, that is, an endomorphism of the tangent bundle satisfying $J^2 = -id$. The complexified tangent bundle $T_{{\mathbb C} } X$ can be decomposed as a direct sum of $T^{(1,0)} X$ and $T^{(0,1)}X$, which are the eigenspaces of $J$ with eigenvalues $\sqrt{-1}$ and $-\sqrt{-1}$, respectively. An almost Hermitian metric $g$ on $X$ is a Riemannian metric that is compatible with $J$, i.e. $g(Y,Z) = g(JY, JZ)$ for any vector fields $Y, Z$. The associated Hermitian form $\omega$ is a real positive definite $(1,1)$-form defined by $\omega(Y, Z) = g(JY, Z)$. For a smooth function $\varphi\in C^\infty(X)$, we define $\ddbar \varphi$ to be 
\bea
\ddbar \varphi = \frac{1}{2} (d(Jd\varphi))^{(1,1)},
\eea 
the $(1,1)$-part of the $2$-form $\frac{1}{2} d(Jd\varphi)$. As in other sections, we write $\omega_\varphi = \omega + \ddbar \varphi$, and let $h_\varphi = \omega^{-1}\cdot \omega_\varphi$ be the endomorphism of the tangent bundle determined by the $(1,1)$-form $\omega_\varphi$ and the metric $\omega$.

\medskip

With this notation, we can consider the same equation on an almost-complex manifold $X$ as in (\ref{eqn:main}), for a given non-linear operator $f:\Gamma\to {\mathbb R_+}$ satisfying the conditions (1-4) in the Introduction, and a given smooth right hand side $e^F$.  We have then the following general $L^\infty$ estimate for fully non-linear equations on almost-complex manifolds:

\begin{theorem}
\label{thm:main5}
Under the above assumptions, and with $\varphi$ normalized to satisfy ${\rm sup}_X\varphi=0$, there exists a constant $C$ depending on $n, \omega, \gamma, p>2n$ and $\| e^{2nF}\|_{L^1 (\log L)^p(X,\omega^n)}$ such that 
$${\mathrm{sup}}_X |\varphi|\le C.$$
\end{theorem}

\bigskip
Only in the particular case when $f$ is the Monge-Amp\`ere operator was an $L^\infty$ bound for $\varphi$ obtained before \cite{CTW}, and an $L^\infty$ bound for $e^{F}$ was required.

\bigskip
{\bf Potentials on almost-K\"ahler manifolds}

Once we leave the realm of K\"ahler geometry, a severe impediment is the lack of a $\partial\bar\partial$-lemma. However, if the failure of the K\"ahler condition is due to the non-integrability of the almost-complex structure, we are in the realm of almost-K\"ahler geometry, and we can still characterize almost-K\"ahler metrics with symplectic form by a potential. 

\medskip
More specifically, let $(X,J)$ be an almost-complex manifold. An almost-K\"ahler metric is a metric $g(U,V)$ which is $J$-compatible, i.e., $g(JU,JV)=g(U,V)$, and whose associated $2$-form $\o_g(U,V)=g(JU,V)$ is closed. Given a tame symplectic form $\Omega$, that is, a closed $2$-form $\Omega$ satisfying $\Omega(U,JU)>0$, we can consider the class of almost-K\"ahler metrics $g$ with $[\o_g]=[\Omega]$. It has been pointed out by Donaldson \cite{D} that there is a natural version of the Calabi-Yau problem in this context, namely find an almost-K\"ahler metric $\tilde g$ with $[\o_{\tilde g}]=[\Omega]$ and given volume form. He showed that a positive answer to this question would have important consequences in symplectic geometry.

\medskip
In \cite{TWY}, it has been shown that the answer is indeed positive if one can establish $L^\infty$ estimates for the potential $\varphi$ defined for any metric $\tilde g$ with $[\o_{\tilde g}]=[\Omega]$ by
\bea
\Delta_{\tilde g}\varphi=2n-2n{\o_{\tilde g}^{n-1}\wedge\Omega\over \o_g^n},
\quad
{\rm sup}_X\varphi=0.
\eea
Such estimates had been reduced first in \cite{W} to an exponential integral estimate for $\varphi$, and then to an $L^p$ estimate in \cite{TWY}. The methods of this paper can readily reduce them to an $L^1$ estimate:

\bigskip
\begin{theorem}\label{thm:main7}
Let $\tilde g$ be an almost-K\"ahler metric on an almost-complex manifold $X$ of dimension $2n$, with $[\o_{\tilde g}]=[\Omega]$ where $\Omega$ is a tame symplectic form. Then there is a constant $C>0$ depending only on the total volume of $\tilde g$ such that
\bea
{\rm sup}_X|\varphi|\leq C(1+\|\varphi\|_{L^1(X, e^{nF}dV_{\tilde g})}),
\eea
where $e^{nF} = \tilde \omega^n / \Omega^n$.
\end{theorem}

This estimate first appeared in our recent survey \cite{GP22} in honor of S.S. Chern. We were informed by V. Tosatti that the methods of \cite{TWY} can also be extended to derive it. Our main interest in including a brief discussion in this paper is to illustrate how a feature of almost-K\"ahler metrics, namely the lack of first derivative terms in the metric in the Laplace-Beltrami operator, is particularly conducive to our methods. 

\bigskip
This paper is a revised and expanded version of a paper that was first posted as arXiv: 2204:12549. We have added new material on gradient terms in $(n-1)$-form equations, on general fully non-linear equations on almost-complex manifolds, and on the Calabi-Yau equation on tame symplectic manifolds, all to show how to adapt our method of local auxiliary Monge-Amp\`ere equations to very different geometric situations.

\section{The key local estimate}\label{section 2 new}
\setcounter{equation}{0}

Fix the Hermitian manifold $(X,\o)$. We assume that there exists $r_0>0$ so that, at any point $x_0\in X$, there exists a local holomorphic coordinate system $z$ centered at $x_0$, with $\o(x_0)=i\sum_{j=1}^n dz^j\wedge d\bar z^j$ and 
\bea\label{eqn:normalization}
{1\over 2}i\p\bar\p |z|^2\leq\o\leq 2 i\p\bar\p |z|^2 \ {\rm in} \ B(x_0,2r_0)
\eea
where $B(x_0,2r_0)$ is the Euclidian ball of radius $2r_0$. When $X$ is compact, with or without boundary, then such an $r_0$ clearly exists, and can be viewed as determined by $X$ and $\o$. Any dependence on $r_0$ can be absorbed into a dependence on $X$ and $\o$.
When $X$ is not compact, the existence of such an $r_0>0$ is an assumption that we always make. Without loss of generality we may assume $r_0\le 1/2$ throughout the paper.

\begin{lemma}
\label{local}\label{thm:local}\label{lemma 1 new}
Let $\varphi$ be a $C^2$ solution on a Hermitian manifold $(X,\o)$ of the equation
\bea
\label{eqn:f}
f(\lambda[h_\varphi])=e^F
\eea
where the operator $f(\lambda)$ satisfies the conditions (1-4) spelled out in Section \S 1. 
Let $x_0$ be any point in $X$, and assume that
\bea
\varphi(z)\geq \varphi(x_0) \ {\rm for}\ z\in B(x_0,2r_0).
\eea
Fix any $p>n$. Then we have
\bea
-\varphi(x_0)\leq C
\eea
where $C$ is a constant depending only on $(X,\o)$, $p$, $\gamma$, $\|e^{nF}\|_{L^1(\log L)^p}$, and $\|\varphi\|_{L^1}$. Here the $L^1$ and $L^1(\log L)^p$ norms are on the ball $B(x_0,2r_0)$ with respect to the volume form $\o^n$.
\end{lemma}

We will assume $-\varphi(x_0)\ge 2$, otherwise Lemma \ref{local} already holds.

\smallskip

We break the proof into several steps. The first step is the most important, and requires a comparison of the solution $\varphi$ of (\eqref{eqn:f}) with the solution of an auxiliary Monge-Amp\`ere equation. More precisely, set
\bea\label{eqn:u s}
u_s(z)=\varphi(z)-\varphi(x_0)+{1\over 2}|z|^2-s
\eea
for each $s$ with $0<s<2r_0^2$. It is convenient to set $\Omega= B(x_0,2r_0)$. Then $u_s>0$ on $\p \Omega$, and thus the set $\Omega_s$ defined by
\bea
\Omega_s=\{z\in \Omega; u_s(z)<0\}
\eea
is an open and nonempty set with compact closure in $\Omega$. Set
\bea
A_s=\int_{\Omega_s}(-u_s)e^{nF}\o^n.
\eea
Formally, the auxiliary equation which we would like to consider is the following Dirichlet problem for the complex Monge-Amp\`ere equation on the ball $\Omega$,
\bea
(i\p\bar\p\psi_s)^n={(-u_s)\chi_{{\mathbb R}_+}(-u_s)\over A_s}e^{nF}\o^n\ {\rm in}\ \Omega,
\quad
\psi_s=0\ {\rm on}\ \p\Omega,
\eea
where $\chi_{{\mathbb R}_+}(x)$ is the characteristic function of the positive real axis, and $\psi_s$ is required to be plurisubharmonic, $i\p\bar\p\psi_s\geq 0$. Since the right hand side would have singularities, we choose a sequence of smooth {\em positive} functions $\tau_k(x): {\mathbb R}\to {\mathbb R}_+$ which converges to the function $x\cdot\chi_{{\mathbb R}_+}(x)$ as $k\to \infty$. Let then the function $\psi_{s,k}$ be defined as the solution of the following Dirichlet problem
\bea
\label{eqn:psi}
(i\p\bar\p \psi_{s,k})^n={\tau_k(-u_s)\over A_{s,k}}e^{nF}\o^n\ {\rm in}\ \Omega, 
\quad 
\psi_{s,k}=0\ {\rm on}\ \p\Omega
\eea
with $i\p\bar\p\psi_{s,k}\geq 0$, and $A_{s,k}$ is defined by
\bea
A_{s,k}=\int_\Omega \tau_k(-u_s) e^{nF} \omega^n.
\eea
By the classic theorem of Caffarelli-Kohn-Nirenberg-Spruck \cite{CKNS}, this Dirichlet problem admits a unique solution $\psi_{s,k}$ which is of class $C^\infty(\bar\Omega)$. The maximum of $\psi_{s,k}$ can only be attained in $\p\Omega$, and thus $\psi_{s,k}\leq 0$. By definition of the constants $A_{s,k}$, we also have $A_{s,k}\to A_s$ as $k\to\infty$, and
\bea
\int_\Omega (i\p\bar\p\psi_{s,k})^n=1.
\eea

\medskip

\begin{lemma}
\label{lm:comparison}
Let $u_s$ be a $C^2$ solution of the fully non-linear equation (\eqref{eqn:f}) and $\psi_{s,k}$ be the solutions of the complex Monge-Amp\`ere equation (\eqref{eqn:psi}) as defined above. Then we have
\bea
\label{comparison}
-u_s\leq \varepsilon (-\psi_{s,k})^{n\over n+1}\ {\rm on}\ \bar\Omega
\eea
where $\varepsilon$ is the constant defined by $\varepsilon^{n+1}=A_{s,k}\gamma^{-1}{(n+1)^n\over n^{2n}}$.
\end{lemma}

\medskip
\noindent
{\it Proof.} We have to show that the function
\bea\label{eqn:test function}
\Phi=  -\varepsilon (-\psi_{s,k})^{n\over n+1}-u_s
\eea
is always $\leq 0$ on $\bar\Omega$. Let $x_{\max}\in \overline{\Omega}$ be a maximum point of $\Phi$. If $x_{\max}\in\overline{ \Omega}\backslash \Omega_s$, clearly $\Phi(x_{\max})\le 0$ by the definition of $\Omega_s$ and the fact that $\psi_{s,k}< 0 $ in $\Omega$. If $x_{\max}\in \Omega_s$, then we have $\ddbar \Phi (x_{\max})\le 0$ by the maximum principle.  

\medskip

Let $G^{i\bar j} = \frac{\partial \log f(\lambda[h])}{\partial h_{ij}} = \frac{1}{f} \frac{\partial f(\lambda[h])}{\partial h_{ij}}$ be the linearization of the operator $\log f(\lambda[h])$. It follows from the structure conditions of $f$ that $G^{i\bar j}$ is positive definite at $h_{\varphi}$ and 
$$\det G^{i\bar j} = \frac {1}{f^n} \det (\frac{\partial f(\lambda[h])}{\partial h_{ij}}) \ge \frac{\gamma}{f(\lambda)^n}.$$ The assumption that $f(\lambda)$ is homogeneous of degree 1 implies that
$$\tr_G \omega_\varphi = \sum_{i,j} G^{i\bar j} (\omega_{\varphi})_{\bar j i} = 1. $$ 
We calculate at the point $x_{\max}$:
\bea
0& \ge& \nonumber (\tr_{G} \ddbar \Phi)(x_{\max}) = \sum_{i,j} G^{i\bar j} \Phi_{\bar j i}\\
& = & \nonumber \frac{n\varepsilon}{n+1} (-\psi_{s,k})^{-\frac{1}{n+1}} (\tr_G \ddbar \psi_{s,k}) + \frac{\varepsilon n}{(n+1)^2} (-\psi_{s,k})^{-\frac{n+2}{n+1}} |\nabla \psi_{s,k}|^2_G \\
 && \nonumber \quad  - \tr_G \omega_\varphi + \tr_G( \omega - \frac 1 2 \ddbar |z|^2  )\\
 & \ge & \nonumber \frac{n\varepsilon}{n+1} (-\psi_{s,k})^{-\frac{1}{n+1}} (\tr_G \ddbar \psi_{s,k}) - 1\\
 & \ge & \nonumber \frac{n^2\varepsilon}{n+1} (-\psi_{s,k})^{-\frac{1}{n+1}} (\det G)^{1/n} [\det ( \omega^{-1}\cdot \ddbar \psi_{s,k})]^{1/n} - 1\\
 & \ge & \nonumber \frac{n^2\varepsilon}{n+1} \gamma^{1/n} (-\psi_{s,k})^{-\frac{1}{n+1}}\frac{(-u_s)^{1/n}}{A_{s,k}^{1/n}} - 1.
 \eea
Here in the fourth line we use \neweqref{eqn:normalization} and in the fifth line we apply the arithmetic geometric inequality. Hence at $x_{\max}$, we have 
$$- u_s\le A_{s,k} \gamma^{-1} \frac{(n+1)^n}{n^{2n} \varepsilon^n} (-\psi_{s,k})^{\frac{n}{n+1}}=  \varepsilon (-\psi_{s,k})^{\frac{n}{n+1}},$$
by the choice of $\varepsilon.$ This implies $\Phi(x_{\max})\le 0$, and the claim is proved. Q.E.D.

\bigskip
Once we have a comparison between $u_s$ and $\psi_{s,k}$ as given in Lemma \ref{lm:comparison}, we no longer need to know the differential equations satisfied by $u_s$ and $\psi_{s,k}$. It is most convenient to summarize then the remaining part of the proof of Theorem 5 in a lemma of general applicability:

\begin{lemma}
\label{general}
Let $0<s<s_0=2r_0^2$. Assume that we have functions $u_s$, $u_s>0$ on $\p\Omega$, and let $\Omega_s$, $A_s$, and $A_{s,k}$ be the corresponding notions as defined by. Assume that
the inequality (\eqref{comparison}) holds, that is,
\bea
-u_s\leq C(n,\gamma)\,A_{s,k}^{1\over n+1}(-\psi_{s,k})^{n\over n+1}\ {\rm on}\ \bar\Omega
\eea
for some constant $C(n,\gamma)$, where $\psi_{s,k}$ are plurisubharmonic functions on $\Omega$, $\psi_{s,k}=0$ on $\p\Omega$, and $\int_\Omega(i\p\bar\p\psi_{s,k})^n=1$. Then
for any $p>n$, we have
\bea
-\varphi(x_0)\leq C(n,\omega, \gamma,p, \|\varphi\|_{L^1(\Omega,\o^n)}).
\eea
\end{lemma}

\bigskip
Clearly, Theorem \ref{local} follows immediately from Lemma \ref{lm:comparison} and Lemma \ref{general}, so we devote the remaining part of this section to the proof of Lemma \ref{general}.

\medskip

The main property of the functions $\psi_{s,k}$ in Lemma \ref{general} that we need is an exponential integral estimate, which can be stated as follows.  Let $\Omega\subset{\mathbb C}^n$ be a bounded pseudoconvex domain. 
Then there are constants $\alpha>0$ and $C>0$ that depend only on $n$ and diam$(\Omega)$ such that 
\begin{equation}
\label{eqn:alpha}
\int_\Omega e^{-\alpha \psi} dV\le C,
\end{equation}
for any $\psi\in PSH(\Omega)\cap C^2(\Omega)\cap  C^0(\overline{\Omega})$ with $\psi|_{\partial \Omega} \equiv 0$, and $\int_\Omega (\ddbar \psi)^n  = 1$.
This exponential estimate is due to Kolodziej \cite{K}, and it is a generalization of the inequality of Brezis-Merle \cite{BM} in complex dimension one. 
While it has some similarity with the $\alpha$-invariant in K\"ahler geometry, it is a very different estimate as it requires information on the Monge-Amp\`ere volume of $\psi$. It was first proved in \cite{K} using pluripotential theory. Recently an elementary and elegant proof was given in \cite{WWZ1} which in particular does not rely on pluripotential theory. 

\bigskip

Returning to the proof of Lemma \ref{general}, we rewrite the comparison between $u_s$ and $\psi_{s,k}$ as
$$\frac{(-u_s)^{(n+1)/n}}{A_{s,k}^{1/n}} \le C(n, \gamma) (-\psi_{s,k}), \mbox{ in } \Omega_s.$$
Take a small $\beta = \beta(n,\alpha, \gamma)>0$ such that $\beta C(n,\gamma)\le \alpha$ where $\alpha>0$ is as in \neweqref{eqn:alpha}. We then get
\begin{equation}
\int_{\Omega_s} \exp\Big( \beta \frac{(-u_s)^{(n+1)/n}}{A_{s,k}^{1/n}} \Big) \omega^n \le \int_\Omega \exp\Big(\beta C(n,\gamma) (-\psi_{s,k}) \Big)\omega^n \le C,
\end{equation}
by the uniform integral estimate \neweqref{eqn:alpha}. We can now take the limit $k\to\infty$ and obtain the following exponential estimate for $u_s$,

\begin{lemma}
For any $0<s<s_0$, the functions $u_s$ satisfy the following inequality
\bea
\label{eqn:m1}
\int_{\Omega_s}{\rm exp}\big\{\beta{(-u_s)^{n+1\over n}\over A_s^{1/ n}}\big\}\omega^n\leq C
\eea
where $\beta=\beta(n,\gamma,r_0)$ and $C=C(n,\gamma,r_0)$ are strictly positive constants depending only on $n$, $\gamma$ and $r_0$.
\end{lemma}

\bigskip
The next step is to show that this exponential inequality for $u_s$ implies the following estimate for $A_s$:

\begin{lemma}
\label{lm:3}
Fix $p>n$, and let $\delta_0={p-n\over pn}>0$. Then we have
\bea
\label{eqn:3}
A_s\leq C_0(\int_{\Omega_s}e^{nF}\o^n)^{1+\delta_0}
\eea
where $C_0$ is a constant depending only on $n, p, \gamma,\omega$ and $\| e^{nF}\|_{L^1(\log L)^p}$.
\end{lemma}

\bigskip
\noindent
{\it Proof.} As in \cite{GPT}, we apply Young's inequality. For our purposes, it is convenient to state Young's inequality as follows. Let $\eta: {\mathbb R}_+\to {\mathbb R}_+$ be a positive strictly increasing function with ${\rm lim}_{x\to 0^+}\eta(x)=0$. Let $\eta^{-1}$ be the inverse function of $\eta$, so that $\eta(\eta^{-1}(y))=y$ and $\eta^{-1}(\eta(x))=x$. Then for any $U,V\geq 0$, we have
\bea
UV\leq \int_0^U\eta(x)dx+\int_0^V \eta^{-1}(y)dy.
\eea
If we make a change of variables $y=\eta(x)$ in the second integral, this inequality can be rewritten as
\bea
UV\leq \int_0^U\eta(x)dx+\int_0^{\eta^{-1}(V)}x\eta'(x)dx
\eea
The first integral on the right hand side can be estimated by $U\eta(U)$, while the second integral can be estimated by
\bea
\int_0^{\eta^{-1}(V)}x\eta'(x)dx\leq \eta^{-1}(V)\int_0^{\eta^{-1}(V)}\eta'(x)dx
=
\eta^{-1}(V)\eta(\eta^{-1}(V))=V\eta^{-1}(V).
\eea
Thus we arrive at the following inequality, which suffices for our purposes,
\bea
UV\leq U\eta(U)+V\eta^{-1}(V).
\eea

We return to the proof of Lemma \ref{lm:3}. We apply the preceding Young's inequality with $\eta(x)=(\log(1+x))^p$, $U=e^{nF(z)}$, and a general positive function $V$. 
Then $\eta^{-1}(y)={\rm exp}(y^{1\over p})-1$, and we find
\bea
\label{Young's}
e^{nF}V\leq e^{nF}(\log(1+e^{nF}))^p+V\big\{{\rm exp}(V^{1\over p})-1\big\}.
\eea
Specializing further to $V=v(z)^p$ for a general non-negative function $v(z)$, we obtain
\bea
\label{Young2}
v(z)^p e^{nF}
\leq e^{nF}(\log(1+e^{nF}))^p+v^p(e^v-1)
\leq  e^{nF}(1+n|F|)^p+C_p e^{2v(z)}.
\eea
Let now 
\bea
v(z)={\beta\over 2}\frac{(-u_s)^{(n+1)/n}}{A_{s}^{1/n}}
\eea
and integrate over $\Omega_s$ to obtain
\begin{equation}\label{eqn:m2}
({\beta\over 2})^p\int_{\Omega_s} \frac{(-u_s)^{(n+1)p/n}}{A_{s}^{p/n}} e^{nF} \omega^n 
\leq \int_{\Omega_s} e^{nF}(1 + n|F|)^p \omega^n + 
\int_{\Omega_s}{\rm exp}\big\{{\beta \frac{(-u_s)^{(n+1)/n}}{A_{s}^{1/n}}}\big\}\o^n.
\nonumber
\end{equation}
By Lemma \ref{lm:3}, the right hand side is bounded by a constant $C(n,p,\beta,\|e^{nF}\|_{L^1(\log L)^p},r_0)$. Thus we have
\begin{equation}\label{eqn:m3}
\int_{\Omega_s} (-u_s)^{\frac{p(n+1)}{n}} e^{n F}\omega^n \le C(n,p,\beta,r_0) A_s^{\frac{p}{n}}.
\end{equation}
This implies the following upper bound for $A_s$: by H\"older's inequality, we can write
\bea
A_s & = & \nonumber \int_{\Omega_s} (-u_s) e^{nF} \omega^n \\
& \le & \nonumber \Big(\int_{\Omega_s} (-u_s)^{p(n+1)/n} e^{nF} \omega^n  \Big)^{n/(p(n+1))} \cdot \Big(\int_{\Omega_s}  e^{nF} \omega^n  \Big)^{1/q} \\
& \le & \nonumber A_s^{1/(n+1)} \cdot \Big(\int_{\Omega_s}  e^{nF} \omega^n  \Big)^{1/q}
\eea
where $q>1$ is chosen such that $\frac{1}{q} + \frac{n}{p(n+1)} = 1$. This is the inequality to be proved in Lemma \ref{lm:3}.  Q.E.D.

\bigskip
We can now derive growth estimates for the function $\phi(s)$ defined by
\bea
\phi(s)=\int_{\Omega_s}e^{nF}\o^n
\eea
following the original strategy of De Giorgi \cite{DG}. We shall actually be interested in estimates as $s\to 0^+$, so the version relevant to our set-up is the one provided by \cite{K}.

\begin{lemma}
\label{lm:4}
The function $\phi(s)$ is a monotone increasing function of $s\in (0,s_0)$ which satisfies ${\rm lim}_{s\to 0^+}\phi(s)=0$, $\phi(s)>0$ for $s\in (0,s_0)$, and
\begin{equation}\label{eqn:m6}
t \phi(s-t) \le C_0( \phi(s))^{1+\delta_0},\quad \mbox{for all } 0< t< s<s_0.
\end{equation}
\end{lemma}

\bigskip
\noindent
{\it Proof.}  The monotonicity of $\phi$ is clear from the definition and the relation $\Omega_{s'}\subset \Omega_s$ if $0<s'< s< s_0$. 
Note that $x_0\in \Omega_s$ for any $s>0$. Being a nonempty open set, $\Omega_s$ has positive measure. So $\phi(s)>0$ for any $s>0$.
Furthermore, $\varphi(z)-\varphi(x_0)+{1\over 2}|z|^2$ is strictly positive in any small punctured neighborhood of $x_0$,
so we have 
$$\cap_{s>0} \overline{ \Omega_s} = \{x_0\}.$$
This shows that $\phi(s)\to 0$ as $s\to 0^+$.

\smallskip

For any $t\in (0,s)$ and any $z\in \Omega_{s-t} \subset \Omega_s$, we have 
$$0> u_{s-t}(z) = u_s(z) + t,\quad \mbox{i.e. } - u_s(z)> t>0.$$
This observation shows that 
\begin{equation}\label{eqn:m5}A_s  = \int_{\Omega_s} (-u_s) e^{nF} \omega^n\ge \int_{\Omega_{s-t}} (-u_s) e^{nF} \omega^n\ge t \int_{\Omega_{s-t}} e^{nF} \omega^n.\end{equation}
Note that $\phi(s) = \int_{\Omega_s}  e^{nF} \omega^n $ for $s\in (0,s_0)$, then equations \neweqref{eqn:3} and \neweqref{eqn:m5} imply that 
\begin{equation}\label{eqn:m6}
t \phi(s-t) \le C_0( \phi(s))^{1+\delta_0},\quad \mbox{for all } 0< t< s<s_0.
\end{equation}
This finishes the proof of Lemma \ref{lm:4}.  Q.E.D.

\begin{lemma}
\label{lm:5}
Let $\phi(s)$ be a function on $(0,s_0)$ satisfying all the properties listed in the previous lemma. There there is a constant $c_0>0$, depending only on $s_0$, $C_0$, and $\delta_0$
so that
\bea
\phi(s_0)\geq c_0.
\eea
\end{lemma}

\bigskip
\noindent
{\it Proof.} We apply the iteration argument in \cite{K}.

\smallskip

Given the initial number $s_0$, we define inductively a {\em decreasing} sequence of positive numbers $\{s_j\}_{j\ge 0}$ as follows: given $s_j$ let
$$s_{j+1} = \sup\{ 0<s< s_j| ~ \phi(s) \le \frac{1}{2}\phi(s_j)  \}. $$
The continuity of $\phi(s)$ implies that $s_{j+1}$ is the largest $s< s_j$ such that $\phi(s) = \frac{1}{2} \phi(s_j)$. Thus we have $\phi(s_{j+1}) = \frac 12 \phi(s_j)$ and $\frac 1 2 \phi(s_j)<\phi(s)\le \phi(s_j)$ for all $s\in (s_{j+1}, s_j)$. Iterating this identity we infer that 
$$\phi(s_j) = 2^{-j} \phi(s_0),\quad \forall j.$$
In particular, $\lim_{j\to\infty} \phi(s_j) = 0$ and Lemma \ref{lm:4} implies that $\lim_{j\to\infty} s_j = 0$, which exists since $\{s_j\}_{j\ge 0}$ is a bounded decreasing sequence of positive numbers. 

\smallskip

From \neweqref{eqn:m6}, we get that for any $s\in (s_{j+1}, s_j)$
$$\frac{1}{2}(s_j - s) \phi(s_j) \le (s_j - s) \phi(s) \le C_0 \phi(s_j)^{1+\delta_0},$$
so $s_j - s\le 2 C_0 \phi(s_j)^{\delta_0}$. Letting $s\to s_{j+1}$ we obtain
$$-s_{j+1} + s_j \le 2 C_0 \phi(s_j)^{\delta_0} \le 2 C_0 2^{-\delta_0 j} \phi(s_0)^{\delta_0}.$$
Taking summation over $j$, we conclude that 
$$s_0 = \sum_{j=0} ^\infty (s_j - s_{j+1}) \le \frac{2C_0}{1-2^{-\delta_0}} \phi(s_0)^{\delta_0}.   $$
Thus we have 
\begin{equation}\label{eqn:lower bound}
\phi(s_0) \ge c_0
\end{equation} for some constant $c_0>0$ depending on $n, \omega, C_0,s_0$. Q.E.D.

\bigskip
\noindent
{\it Proof of Lemma 1}. 
Note that on $\Omega_{s_0}$, $u_{s_0} = \varphi + \frac 1 2 |z|^2 - \varphi(x_0) - s_0<0$, that is,
\begin{equation}\label{eqn:f1}
-\varphi - \frac 1 2 |z|^2 > -\varphi(x_0) - s_0 >1 \ \quad \mbox{on }\Omega_{s_0},
\end{equation}
and the last inequality holds from our assumption that $s_0< 1/2$ and $\varphi(x_0)< -2$. From \neweqref{eqn:f1} we deduce that
\begin{equation}\label{eqn:f2}
\log \frac{  -\varphi - 2^{-1} |z|^2  }{( -\varphi(x_0) - s_0  )^{1/2}} > \log ( -\varphi(x_0) - s_0  )^{1/2}>0\ \quad \mbox{on }\Omega_{s_0}.
\end{equation}
Integrating \neweqref{eqn:f2} against the measure $e^{nF} \omega^n$ over $\Omega_{s_0}$, we deduce that
\bea
\label{eqn:f3}
 &&\nonumber \log ( -\varphi(x_0) - s_0  )^{1/2}\cdot \phi(s_0)  = \log ( -\varphi(x_0) - s_0  )^{1/2} \int_{\Omega_{s_0}} e^{nF} \omega^n\\
 & < &\nonumber  \int_{\Omega_{s_0}}  \log \frac{  -\varphi - 2^{-1} |z|^2  }{( -\varphi(x_0) - s_0  )^{1/2}}  e^{nF} \omega^n
 \eea
 We now apply the same Young's inequality (\ref{Young's}) as before, but this time with the following choice of the factor $V$,
 \bea
 V=\log \frac{  -\varphi - 2^{-1} |z|^2  }{( -\varphi(x_0) - s_0  )^{1/2}}=:\log\,W.
 \eea
 The inequality becomes
 \bea
( \log\,W)
 \,e^{nF}
 \leq
 e^{nF}(\log(1+e^{nF}))^p+(\log\,W)\big\{{\rm exp}(\log\,W)^{1\over p}-1\big\}
 \eea
 Fix now $p>1$. A rough upper bound of the two terms on the right hand side which suffices for our purposes is the following
 \bea
 e^{nF}(\log(1+e^{nF}))^p\leq e^{nF}(1+|nF|)^p
 \eea
 and
 \bea
 (\log\,W)\big\{{\rm exp}(\log\,W)^{1\over p}-1\big\}
& \leq& (\log\,W)\big\{{\rm exp}({1\over p}\log\,W+C_p)-1\big\}\nonumber\\
&\leq & C_p'(\log\,W)W^{1\over p}
\leq C_p''(1+|W|).
\eea
Substituting back this inequality in the initial inequality (\ref{eqn:f3}), and estimating the integral over $\Omega_s$ by the full integral over $\Omega$,
we find
\bea
\nonumber \log ( -\varphi(x_0) - s_0  )^{1/2} \cdot\phi(s_0)
&\leq&
\|e^{nF}\|_{L^1(\log L)^p}
+C_p''\int_\Omega (1+{-\varphi(z)-{1\over 2}|z|^2\over (-\varphi(x_0)-s_0)^{1\over 2}})\o^n
\nonumber\\
&\leq &
\|e^{nF}\|_{L^1(\log L)^p}+C_p''(C_1(\o)+{\|\varphi\|_{L^1(\Omega)}+C_2(\omega)\over (-\varphi(x_0)-s_0)^{1\over 2}})
\nonumber
\eea
In particular,
\bea
 \log ( -\varphi(x_0) - s_0  )^{1/2}\cdot \phi(s_0)
 \leq {1\over 2}{\max}\big\{\|e^{nF}\|_{L^1(\log L)^p}+C_p''C_1(\o),C_p''{\|\varphi\|_{L^1(\Omega)}+C_2(\o)\over  (-\varphi(x_0)-s_0)^{1\over 2}}\big\}.
 \nonumber
 \eea
Since $\phi(s_0)$ is bounded from below by Lemma \ref{lm:5}, it follows immediately that $-\varphi(x_0)-s_0$ is bounded above by a constant depending only on $X,\o,n,p, \|e^{nF}\|_{L^1(\log L)^p}$ and $\|\varphi\|_{L^1}$. The proof of Theorem \ref{local} is complete. Q.E.D.
 
\medskip

We remark that compared to the compact K\"ahler case in \cite{GPT}, the constant $C>0$ in Theorem \ref{thm:local} does not depend on the ``energy'' $E(\varphi) = \int_X(-\varphi) e^{nF} \omega^n$. This is because the function $\psi_{s,k}$ in the auxiliary equation \neweqref{eqn:psi} is {\em local} and {\em plurisubharmonic} in $\Omega_s$, while its counterpart in \cite{GPT} is {\em global} and $\omega$-plurisubharmonic. Thus we do not need to shift the function $-\psi_{s,k}$ in \neweqref{eqn:test function} by a constant depending on $E(\varphi)$. The proof of Theorem \ref{thm:local} is in particular simpler since it does not require an argument for {\em bounding the energy by the entropy} as in \cite{GPT} .

\medskip
We also observe one feature of the auxiliary equation which sets it apart from many other auxiliary equations in the literature such as in \cite{CC}: it depends itself on the unknown function $u$. This may seem to limit its applicability, but as we see from Lemma \ref{lemma 1 new}, it turns out that for our purposes, very little is needed from the auxiliary functions $\psi_{s,k}$.

\section{The case of Hermitian compact manifolds without boundary}\label{section 3 new}
\setcounter{equation}{0}

Once we have Lemma \ref{local}, it is easy to establish all the theorems listed in the Introduction. In this section, we begin with the case of $X$ Hermitian compact without boundary, and give the proof of Theorem \ref{thm:main1}.

\bigskip
\noindent
{\it Proof of Theorem \ref{thm:main1}.} Since $X$ is compact without boundary, the function $\varphi(z)$ must attain its minimum at some interior point $x_0$. By Lemma \ref{local}, $-\varphi(x_0)$ is bounded above by a constant depending only on $X,\o,p,n, \gamma, \|e^{nF}\|_{L^1(\log L)^p}$, and $\|\varphi\|_{L^1(B(x_0,2r_0))}$. Thus the proof is complete once we have proved the following lemma giving a bound on the  $L^1$-norm of $\varphi$ for which $\lambda[h_\varphi]\in \Gamma\subset \Gamma_1$.

\begin{lemma}\label{lemma 3}
For any $\varphi\in C^2(X)$ such that $\lambda[h_\varphi]\in \Gamma\subset \Gamma_1$ and $\sup_X \varphi =0$, there exists a uniform constant $C>0$ depending only $n, \omega$ such that 
$$\int_X (-\varphi) \omega^n \le C.$$
\end{lemma}
\noindent {\em Proof.} Note that $\lambda[h_\varphi]\in \Gamma_1$ means that $\tr_\omega (\omega+ \ddbar \varphi)>0$, that is, in local coordinates, the following holds
\begin{equation}\label{eqn:n1}g^{i\bar j} \frac{\partial^2 (-\varphi)}{\partial z_i\partial \bar z_j} < n.\end{equation}
If $\omega$ is K\"ahler, the lemma follows easily from the Green's formula. In general, we will apply local elliptic estimates to $-\varphi$ over a finite cover of $(X,\omega)$. More precisely, we fix three sets of open covers $\{U_\ell\}_{\ell =1}^N$, $\{U_\ell '\}_{\ell = 1}^N$ and $\{U_\ell ''\}_{\ell = 1}^N$, such that $U_\ell\subset \subset U'_\ell\subset \subset U_\ell''$ for each $\ell$. We can also view $U_\ell$ as a Euclidean ball of radius $r_\ell$ and $U_\ell', U_\ell''$ as concentric balls with radii $2 r_\ell, 3 r_\ell$, respectively. The choice of these covers is fixed and depends only on $(X,\omega)$. Assume $\sup_X \varphi=0 = \varphi(p_1)$ is achieved at some point $p_1\in U_1\subset U_1'$.

We can apply the standard elliptic estimate (see e.g. Theorem 8.18 in \cite{GT}) to the equation of $(-\varphi)$ as in \neweqref{eqn:n1} on $U_1'\subset U_1''$, to conclude 
\begin{equation}\label{eqn:n2}
\int_{U_1'} (-\varphi) \omega^n \le C [ (-\varphi)(p_1) + 1   ] = C_1,
\end{equation}
for some constant $C_1>0$ depending on $U_1'\subset U_1''$ and $\omega$ only. If $U_1'\cap U_2' \neq \emptyset$, then from \neweqref{eqn:n2}  we get
$$
\int_{U_1'\cap U_2'} (-\varphi) \omega^n \le C_1,
$$
which by mean-value theorem implies that there is some point $p_2\in U_1'\cap U_2'$ such that $-\varphi(p_2)\le C_1'$ for some $C_1'>0$. We can now apply again the elliptic estimate to $(-\varphi)$ which satisfies \neweqref{eqn:n1} in $U_2'\subset U_2''$ to conclude that
$$\int_{U_2'} (-\varphi) \omega^n \le C[ (-\varphi)(p_2) + 1  ]\le C_2,$$
for some uniform constant $C_2>0$. We can repeat this argument finitely many times, to deduce that $\int_{U_\ell'} (-\varphi) \omega^n \le C_\ell$ for each $\ell$, since $\{U_\ell'\}_{\ell = 1}^N$ is an open cover of $(X,\omega)$, which is assumed to be connected. From this, it is easy to see the bound of $L^1$-norm of $(-\varphi)$. Lemma \ref{lemma 3} is proved. 

\bigskip
The proof of Lemma \ref{local} already contains a proof of Corollary \ref{cor:main}. More precisely, $\varphi$ is by definition a solution on $X$ of the complex Monge-Amp\`ere equation
$$f(\lambda[h_\varphi]) = \Big(\frac{\omega_\varphi^n}{\omega^n}\Big)^{1/n} = e^F$$
which does satisfy the conditions (1-4) described in the Introduction. Let $x_0$ be a point in $X$ where $\varphi$ attains its minimum. Then Lemmas \ref{lm:comparison}-\ref{lm:4} apply, implying 
that there exists a constant $c_0>0$ depending only on $n, p, \omega, $ and $K$ such that 
$$\phi(s_0) = \int_{\Omega_{s_0}} e^{nF} \omega^n \ge c_0.$$
It is clear that $\phi(s_0)\le \int_X (\omega+\ddbar \varphi)^n$. Q.E.D.

\bigskip

We now provide a proof of the inequality \neweqref{eqn:vol lower}, without the extra assumption on the $p$-th entropy of $e^{nF}$, but under an additional condition that $X$ admits a smooth {\em closed} real $(1,1)$-form $\beta$ whose Bott-Chern class is nef and big. The proof of Lemma \ref{lemma Demailly} is motivated by an observation of Tosatti \cite{T} using the Morse inequality of Demailly \cite{D}.

\begin{lemma}\label{lemma Demailly}
Under the assumption above, there exists a constant $c>0$ that depends on $n, \omega$ and $\beta$ such that 
$$\int_X (\omega + \ddbar \varphi)^n \ge c,$$ for any $C^2$-function $\varphi\in PSH(X,\omega)$.
\end{lemma}

\noindent{\em Proof of Lemma \ref{lemma Demailly}. } By a weak transcendental Morse inequality of Demailly \cite{D}, one can show (see e.g. Theorem 3.1 in \cite{T}) that for {\em any} $C^2$ function $u$, the following holds, 
\begin{equation}\label{eqn:Demailly}
0< \int_X \beta^n = \int_X \beta_u^n \le \int_{X(\beta_u, 0)} \beta_u^n,
\end{equation}
where $\beta_u = \beta + \ddbar u$, and $X(\beta_u, 0)$ denotes the subset of $X$ where the $(1,1)$-form $\beta_u$ is nonnegative definite.

\smallskip

Since $\beta$ is a fixed $(1,1)$-form and $\omega$ is positive definite, there exists a constant $\Lambda>0$ such that $  \beta \le \Lambda \omega$ on $X$. Applying \neweqref{eqn:Demailly} to $u = \Lambda \varphi$ for any $\varphi \in C^2\cap PSH(X,\omega)$ gives
\bea
0<\int_X \beta^n &\le &\int_{X(\beta_u, 0)} (\beta + \ddbar \Lambda \varphi)^n
\le \int_{X(\beta_u, 0)} (\Lambda \omega + \Lambda \ddbar \varphi )^n\nonumber\\
&\le & \Lambda^n \int_X (\omega+ \ddbar \varphi)^n,\nonumber
\eea
since $0\le \beta + \ddbar \Lambda \varphi \le \Lambda \omega + \ddbar \Lambda \varphi$ over the set $X(\beta_u, 0)$.  Q.E.D.

\section{The case of Hermitian manifolds with boundary}
\setcounter{equation}{0}

We consider now the case of the Dirichlet problem for the equation (\ref{eqn:DMA}) for an operator $f(\lambda)$ satisfying the conditions (1-4), and give the proof of Theorem \ref{thm:main2}.

\bigskip
\noindent
{\it Proof of Theorem \ref{thm:main2}}. In the case of the Dirichlet problem, we cannot normalize the solution to satisfy ${\sup}_X\varphi=0$, and we need to establish an upper bound for $\varphi$. However, this can be readily done by comparing $\varphi$ to a solution of the Dirichlet problem for the Laplace-Beltrami equation on $X$. More precisely, let $h$ be the solution of the equation
\bea
\label{LB}
\Delta_\o h=n+1\ {\rm in}\ X,
\quad h=0\ \ {\rm on}\ \p X.
\eea
where $\Delta_\o h =n(i\p\bar\p h\wedge\o^{n-1})\o^{-n}$ is the complex Laplacian defined by the Hermitian metric $\o$. In the K\"ahler case, $\Delta_\o$ coincides with the Laplace-Beltrami equation, and the unique solvability of the corresponding Dirichlet problem is well-known. It has been shown by Z. Lu \cite{Lu} that the unique solvability of the Dirichlet problem still holds for general Hermitian metrics $\o$ on compact manifolds $X$ with smooth boundary. Thus the function $h$ exists and is unique, and a dependence on it can be absorbed into a dependence on $(X,\o)$.

\medskip
Now $\lambda[h_\varphi]\in\Gamma\subset\Gamma_1$ implies that $\Delta_\o\varphi\geq -n$. Thus $\Delta(\varphi+h)\geq 1$ on $X$, which implies that $\varphi+h$ must attain its maximum on the boundary. Hence
\bea
{\sup}_X(\varphi+h)
\leq {\sup}_{\p X}(\varphi+h)
={\rm sup}_{\p X}\rho.
\eea
This gives an upper bound for $\varphi$ by a constant depending only on $X,\o$, and $\rho$.

\medskip
We can now replace $\varphi$ by $\varphi-{\sup}_X\varphi$, and thus $\varphi$ is a non-positive solution
to \neweqref{eqn:DMA}. Suppose $\inf_X \varphi = \varphi(x_0)$ for some $x_0\in \overline X$. If $x_0\in \partial X$, we are done. Otherwise we assume $x_0\in X$ and 
\begin{equation}\label{eqn:new assumption}\varphi(x_0) \le \inf_{\partial X} \rho - 1.\end{equation} Take a local coordinates system $(B(x_0, 2r_0), z)$ centered at $x_0$ and \neweqref{eqn:normalization} holds. 

\smallskip

If $B(x_0, 2r_0)\subset\subset X$, we can apply directly Theorem \ref{local} and Lemma \ref{lm:4}, which gives the desired bound
for $-\inf_X\varphi$, modulo an $L^1(X, \omega)$-bound of $(-\varphi)$. 

\smallskip

If $B(x_0, 2r_0)\cap \partial X \neq \emptyset$, we let $s_0 = 2 r_0^2$, $s\in (0, s_0)$, and consider the function 
\begin{equation}
u_s(z) = \varphi(z) - \varphi(x_0) + \frac 12 |z|^2 - s, \quad \forall z\in \Omega = B(x_0, 2r_0)\cap X.
\end{equation}
We observe that on $\partial B(x_0, 2r_0) \cap X$, $u_s\ge 2 r_0^2 - s>0$; and on the other half of $\partial \Omega$, $\forall z \in B(x_0, 2r_0)\cap \partial X$, we have by \neweqref{eqn:new assumption}
$$u_s(z) = \rho(z) - \varphi(x_0) + \frac 12 |z|^2 - s\ge 1 - s>0.$$
Thus we have $u_s|_{\partial \Omega}>0$. We can extend the function $u_s$ to a smooth function on $B(x_0, 2r_0)$ satisfying $u_s>0$ on $B(x_0,2r_0)\backslash \Omega$. Since only the region when $u_s<0$ is under consideration, the estimates in \S 2 are independent of the extension of $u_s$. We also extend $F$ to be smooth function on $B(x_0,2r_0)$, and it will be clear from \neweqref{eqn:aux new} that the extension is again irrelevant. 

\smallskip

Denote the sub-level set of $u_s$, $\Omega_s = \{z\in \Omega| u_s(z)<0\}$. By the discussions above, $\Omega_s\subset \subset X$. As in Section \S 2, we consider the following (auxiliary) Dirichlet boundary problem (with the extended functions $u_s$ and $F$ above)
\begin{equation}\label{eqn:aux new}
(\ddbar \psi_{s,k})^n = \frac{\tau_k(-u_s)}{A_{s,k}} e^{nF} \omega^n,\quad \mbox{in }B(x_0, 2r_0),
\end{equation}
and $\psi_{s,k}\equiv 0$ on $\partial B(x_0, 2r_0)$. Here $\tau_k$ and $A_{s,k}$ are the same as in Section \S 2. We note that the calculations in Section \S 2 are performed essentially in $\Omega_s$, which is strictly contained in $X$. Therefore, we still obtain Theorem \ref{local} and Lemma \ref{lm:4}, and hence a lower bound of $\varphi(x_0) = \inf_X \varphi$, depending on an $L^1$ bound of $(-\varphi)$. We remark that because $\partial X\neq \emptyset$, the {\em a priori} $L^1$-norm of $\varphi$ as in Lemma \ref{lemma 3} in general does not hold.   

\bigskip

We now focus on complex Monge-Amp\`ere and Hessian equations with $\omega$ being K\"ahler, i.e. $d\omega =0$. We rewrite the equation \neweqref{eqn:main} as ($1\le k\le n$)
\begin{equation}\label{eqn:4.6}
(\omega+\ddbar \varphi)^k\wedge \omega^{n-k} = e^{\tilde F} \omega^n \mbox{ in }X,\,\, \varphi = \rho \mbox{ on }\partial X.
\end{equation}
Here we denote $\tilde F = k F$ to simplify notations. We will write $\omega_\varphi = \omega + \ddbar \varphi\in \Gamma_k$. When $k = n$, \neweqref{eqn:4.6} reduces to the Monge-Amp\`ere equation. Since $\varphi$ is always bounded above depending on $\omega$ and $ \rho$, we may assume $\sup_X \varphi \le 0$, and that $\inf_X \varphi = \varphi(x_0) \le \inf_{\partial X}\rho - 1$ for some $x_0\in X^\circ$, otherwise we are done. 

\smallskip

 For simplicity of notations, we consider the function $\tilde \varphi = \min(0, \varphi - \inf_{\partial X} \rho)$, which is Lipschitz in $X$ and $\tilde \varphi|_{\partial X} \equiv 0$. We calculate
 \bea
\nonumber
&&\int_X \log (-\tilde \varphi + 1) ( \omega_\varphi^k - \omega^k )\wedge \omega^{n-k}\\
& = \nonumber& \int_X \log(-\tilde \varphi + 1) \ddbar \varphi\wedge ( \omega_\varphi^{k-1} + \omega_\varphi^{k-2}\wedge \omega + \cdots + \omega^{k-1}  ) \wedge \omega^{n-k}\\
& = &\nonumber \int_X \frac{1}{-\tilde \varphi + 1} i\partial \tilde \varphi\wedge \bar \partial \tilde \varphi \wedge ( \omega_\varphi^{k-1} + \omega_\varphi^{k-2}\wedge \omega + \cdots + \omega^{k-1}  ) \wedge \omega^{n-k}\\
&& \label{eqn:4.7} + \int_{\partial X} \log(-\tilde \varphi+ 1) \bar \partial \varphi \wedge ( \omega_\varphi^{k-1} + \omega_\varphi^{k-2}\wedge \omega + \cdots + \omega^{k-1}  ) \wedge \omega^{n-k} \\
&&\label{eqn:4.8} - \int_X \log (-\tilde \varphi + 1) \bar \partial \tilde \varphi \wedge \partial \Big(( \omega_\varphi^{k-1} + \omega_\varphi^{k-2}\wedge \omega + \cdots + \omega^{k-1}  ) \wedge \omega^{n-k} \Big)\\
&\ge & \frac{4}{n} \int_X\Big |\nabla \sqrt{-\tilde \varphi + 1}\Big |^2_\omega \omega^n\nonumber
\eea
where the integral in \neweqref{eqn:4.7} vanishes since $\log(-\tilde \varphi + 1)|_{\partial X} \equiv 0$, and that in \neweqref{eqn:4.8} is zero since we assume $\omega$ is K\"ahler.

\smallskip

On the other hand, the left-hand-side of the equation above satisfies
\bea\nonumber 
 \int_X \log (-\tilde \varphi + 1) ( \omega_\varphi^k - \omega^k )\wedge \omega^{n-k} = \int_X \log(-\tilde \varphi + 1) (e^{\tilde F}-1)\omega^n
 \le\int_X \log(-\tilde \varphi + 1) e^{\tilde F}\omega^n
\eea
since $\log(-\tilde\varphi+1)\geq 0$. We can now apply Young's inequality as in (\ref{Young's}) with $p=1$, $v={1\over 4}\log(-\tilde\varphi+1)$ to write
\bea
\nonumber
\log(-\tilde \varphi + 1) e^{\tilde F}
\leq 
4 \,e^{\tilde F} ( 1+ |\tilde F|) +C\sqrt{-\tilde\varphi+1}
\eea
and hence
\bea
\nonumber
 \int_X \log (-\tilde \varphi + 1) ( \omega_\varphi^k - \omega^k )\wedge \omega^{n-k}
& \leq&
4 \int_X e^{\tilde F} ( 1+ |\tilde F|) \omega^n + C \int_X\sqrt{-\tilde\varphi+1}\omega^n\nonumber\\
&\leq& C + C \int_{X} \sqrt{-\tilde \varphi + 1} \omega^n.\nonumber
\eea
Combining the above, we obtain 
$$\int_{X} \Big|\nabla ( \sqrt{-\tilde \varphi + 1} - 1  ) \Big|^2_\omega\omega^n \le C + C \int_X ( \sqrt{-\tilde \varphi + 1} - 1  ) \omega^n.$$
Observing that the function $(\sqrt{-\tilde\varphi + 1} - 1) \equiv 0$ on $\partial X$, we can apply the Poincar\'e inequality on $({\overline X}, \omega)$ to get ($C_P$ below is the Poincar\'e constant)
\bea\nonumber \int_X ( \sqrt{-\tilde \varphi + 1} - 1  )^2 \omega^n & \le & C_P \int_X \Big|\nabla ( \sqrt{-\tilde \varphi + 1} - 1) \Big|_\omega^2\omega^n \le C + C\int_X ( \sqrt{-\tilde \varphi + 1} - 1) \omega^n\\
&\le &C + C\Big(\int_X ( \sqrt{-\tilde \varphi + 1} - 1)^2 \omega^n\Big)^{1/2}\nonumber,
\eea
from which we easily get $\int_X ( \sqrt{-\tilde \varphi + 1} - 1)^2\omega^n\le C$, hence $$\int_X (-\tilde \varphi)\omega^n \le C.$$ Since $\tilde \varphi$ differs from $\varphi$ by a uniform constant that depends only on $\inf_{\partial X} \rho$, this proves the $L^1(X,\omega^n)$-bound of the function $|\varphi|$.  Q.E.D.

\section{The case of non-compact Hermitian manifolds}
\setcounter{equation}{0}

We now turn to the proof of Theorem \ref{thm:main3}, which is analogous to that of Theorem \ref{thm:main2}, so we only point out the very minor differences.

\medskip
\noindent
{\it Proof of Theorem \ref{thm:main3}}. First note we are allowing a dependence on the (local) $L^1$ norm of $\varphi$ in the constant $C$ in Theorem \ref{thm:main3}. With such an allowance, the upper bound of $\varphi$ is local since $\Delta_\omega \varphi \ge -n$ and it follows from the standard elliptic estimates (e.g. Theorem 8.18 in \cite{GT}). 

Next we derive the lower bound for $\varphi$. Let ${\rm liminf}_{z\to\infty}\varphi> -B$ for some positive constant $B$ which can be viewed as part of the data. Then there is a compact set $K\subset\subset X$ with $\varphi\geq -B-1$ for $z\in X\setminus K$. If $\varphi\geq -B$ on $K$, we are done. Otherwise, the infimum of $\varphi$ over $X$ must take place at some interior point $x_0$ of $K$.
We can now apply Theorem \ref{local}, and obtain the desired bound for $-\varphi(x_0)=-{\inf}_X\varphi$. Q.E.D.

\section{The case of $(n-1)$-form type nonlinear equations}
\setcounter{equation}{0}

We come now to the $(n-1)$ form type nonlinear equation and the proof of Theorem \ref{thm:main4}. 

\medskip
\noindent
{\it Proof of Theorem \ref{thm:main4}.} We first fix a constant $C_\theta = C_\theta(\omega_h, \theta)$ such that $i\theta\wedge \bar \theta \le C^2_\theta \omega_h$.
Recall that $\varphi$ is a $C^2$ solution to the $(n-1)$-form nonlinear equation \neweqref{eqn:HLMA}. For notational convenience, we denote
$$\tilde \omega = \omega_h + \frac{1}{n-1} ( (\Delta_\omega \varphi) \omega - \ddbar \varphi ) + \chi[\varphi]$$
to be the $(1,1)$-form on the left-side of \neweqref{eqn:HLMA}, and $\tilde h_\varphi = \omega^{-1}\cdot \tilde \omega$. Write $G^{i\bar j} = \frac{\partial \log f(\lambda[\tilde h_\varphi])}{\partial \tilde h_{i\bar j}}$ to be linearized operator of $\log f(\lambda[\cdot])$ which is positive definite by (3) of the assumptions on $f$, and we will denote $\tr_G \alpha = G^{i\bar j} \alpha_{i\bar j}$ for any $(1,1)$-form $\alpha$. 
Write $g = (g_{i\bar j})$  to be the metric associated to the background Hermitian metric $\omega$ and $g^{i\bar j}$ its inverse. 

\smallskip

As in \cite{TW1, KLW}, we define a tensor $\Theta^{i\bar j}$ by 
$$\Theta^{i\bar j} = \frac{1}{n-1} \Big( (\tr_{G} \omega) g^{i\bar j} - G^{i\bar j}   \Big).$$
\begin{lemma}\label{lemma 3.1}
The tensor $\Theta^{i\bar j}$ is positive definite at any point, and satisfies
\begin{equation}\label{eqn:3.1}
\det(\Theta^{i\bar j}) \ge \frac{\gamma}{e^{nF}} \det(g^{i\bar j}), \mbox{\, on }X.
\end{equation}
\end{lemma}
\noindent{\em Proof. } Given any $x_0\in X$, we can choose local holomorphic coordinates around $x_0$ such that $g_{i\bar j}|_{x_0} = \delta_{ij} $ and $\tilde \omega|_{x_0} = \lambda_i \delta_{ij}$ is diagonal. Note that $\lambda = (\lambda_1,\ldots,\lambda_n)\in \Gamma$ and $G^{i\bar j}|_{x_0} = \frac{1}{f} \frac{\partial f(\lambda)}{\partial \lambda_i} \delta_{ij} = \mu_i \delta_{ij}$ for $\mu_i = \frac{1}{f} \frac{\partial f(\lambda)}{\partial \lambda_i}>0$. 
Therefore, at $x_0$, we have
$$\Theta^{i\bar j} = \frac{1}{n-1} (\sum_{k\neq i}\mu_k) \delta_{ij}$$ which is clearly positive definite. Furthermore, we have
\bea\det (\Theta^{i\bar j})|_{x_0} &= &\nonumber  \frac{1}{(n-1)^n} \prod_{i=1}^n \sum_{k\neq i} {\mu_k}\ge \prod_{i=1}^n \Big( \prod_{k\neq i} {\mu_k}   \Big)^{1/(n-1)}\\
 &= & \nonumber\prod_{i=1}^n {\mu_i} = \frac{1}{f^n} \prod_{i=1}^n \frac{\partial f}{\partial \lambda_i} \ge \frac{\gamma}{f^n},
\eea
where the inequality follows from the arithmetic-geometric (AG) inequality and the last equality follows from the structural condition \neweqref{eqn:structure} on $f$. The proof of the lemma is complete since $x_0\in X$ is an arbitrary point in $X$.
Q.E.D.

\medskip

Next, we take a minimum point $x_0$ of $\varphi$, i.e. $\varphi(x_0) = \inf_X\varphi$. Let $(z_1,\cdots,z_n)$ be a local holomorphic coordinate system centered at $x_0$, such that on the Euclidean ball $\Omega = B(x_0, 2 r_0)\subset {\mathbb C}^n$ defined by these coordinates, the following holds
\begin{equation}\label{eqn:3.2}
\frac{1}{2} \ddbar |z|^2 \le \omega \le 2 \ddbar |z|^2.
\end{equation}
We fix a small constant $\epsilon'>0$ which will be determined later. 
Let $z_i = x_{2i-1} + \sqrt{-1} x_{2i}$. Hence $\{x_a\}_{a=1}^{2n}$ is a local real coordinates system. 
Let $s_0 =  \epsilon' r_0^2$. For any $s\in (0,s_0)$, we define a function $u_s$ similar to that in \neweqref{eqn:u s}.
\begin{equation}\label{eqn:u s 1}
u_s(z) : = \varphi(z) - \varphi(x_0) + \epsilon' |z|^2 - s,\quad \forall z\in \Omega := B(x_0, 2r_0).
\end{equation}
Note that when $z\in  \Omega\backslash B(x_0, r_0)$, $u_s(z)\ge \epsilon' r_0^2 - s>0$. Therefore the sub-level set of $u_s$,
$\Omega_s : = \{z\in \Omega|~ u_s(z)<0\}$
is relatively compact in $\Omega$ and is contained in $B(x_0, r_0)$. By definition $\Omega_s$ is a nonempty open set, since $x_0\in \Omega_s$. We consider the following real Monge-Amp\`ere equation
$$\det \Big(\frac{\partial^2 \psi_{s,k}}{\partial x_a \partial x_b}\Big) = \frac{\tau_k(-u_s)}{A_{s,k}} e^{2nF} (\det g_{i\bar j})^2, \mbox{ in }\Omega,$$ 
and $\psi_{s,k}$ is convex such that $\psi_{s,k} = 0$ on $\partial \Omega$. Here $A_{s,k}>0$ is chosen so that $\int_{\Omega} \det\Big( \frac{\partial^2 \psi_{s,k}}{\partial x_i \partial x_j}\Big)  dx = 1$, i.e. $A_{s,k} = \int_{B(x_0, 2r_0)} \tau_k(-u_s) e^{2nF} (\det g_{i\bar j})^2 dx \to \int_{\Omega_s} (-u_s) e^{2nF} (\det g_{i\bar j})^2 dx=: A_s.$ Note that $A_{s}\le C(\omega) s_0 \| e^{2nF}\|_{L^1(X,\omega^n)}$, so, we may assume $A_{s,k}$ is uniformly bounded, by taking $k$ large enough.  
Denote $\beta_n$ the volume of the unit ball in ${\mathbb R}^{2n}$.
\begin{lemma}\label{lemma 11}
There exists a constant $C_0=C_0(n)>0$ such that 
\begin{equation}\label{eqn:Feb18}-{\mathrm{inf}}_{B(x_0, 2r_0)} \psi_{s,k} \le C_0 r_0,\quad {\rm sup}_{B(x_0, r_0)} |\nabla \psi_{s,k}|\le C_0.\end{equation}
\end{lemma}
\noindent {\em Proof.} By the definition of $\psi_{s,k}$, $ \int_{B(x_0, 2r_0)} \det \Big( \frac{\partial^2 \psi_{s,k}}{\partial x^i \partial x^j}\Big) dx = 1$, it follows from the standard ABP maximum principle (see e.g. \cite{GT})
$$-\inf_{B(x_0, 2r_0)} \psi_{s,k} \le -\inf_{\partial B(x_0, 2r_0)} \psi_{s,k} + \frac{4 r_0}{\beta_n} \Big[ \int_{B(x_0, 2r_0)}\det  \Big( \frac{\partial^2 \psi_{s,k}}{\partial x_i \partial x_j}\Big) dx \Big]^{1/2n} = \frac{4 r_0}{\beta_n} .$$
\noindent To see the second inequality, for any fixed point $x\in B(x_0, r_0)$, denote $Y = \frac{D\psi_{s,k}(x)}{|D\psi_{s,k}(x)|}$ to be the unit vector in the direction of $D\psi_{s,k}(x)$ (if $D\psi_{s,k}(x) = 0$, there is nothing to prove, so here we assume $D\psi_{s,k}(x)\neq 0$). Then clearly $|D\psi_{s,k}(x)| = D\psi_{s,k}(x) \cdot Y$. Consider the half line $L:$ $0\le t\mapsto x + t Y$ which intersects $\partial B(x_0, r_0)$ and $\partial B(x_0, 2r_0)$ at $L(t_1)$, $L(t_2)$, respectively. We have $t_2 - t_1\ge r_0$ and $0>\psi_{s,k}(L(t_1))\ge -\frac{4}{\beta_n} r_0$ and $\psi_{s,k}(L(t_2)) = 0$. Then by the convexity of the function $t\mapsto \psi_{s,k}(L(t))$, we have 
$$|D\psi_{s,k}(x)| = \psi_{s,k}(L(t))'|_{t= 0}\le \frac{\psi_{s,k}(L(t_2)) - \psi_{s,k}(L(t_1))}{t_2 - t_1} \le \frac{4}{\beta_n}.$$
Taking supremum over all $x\in B(x_0, r_0) $ finishes the proof of the lemma. Q.E.D.

\smallskip

It holds that 
\begin{equation}\label{eqn:Blocki}\det (\psi_{s,k})_{i\bar j} \ge 2^{-n} \Big [ \det \Big(\frac{\partial^2 \psi_{s,k}}{\partial x_a \partial x_b}\Big)  \Big]^{1/2} = 2^{-n} \frac{\tau_k(-u_s)^{1/2}}{A_{s,k}^{1/2}} e^{nF} \det (g_{i\bar j}).\end{equation}

The key is now to establish the analogue of Lemma \ref{lm:comparison}, but where $u_s$ is now constructed from the above form-type nonlinear equation, instead of solutions of the non-linear scalar equation $f(\lambda)$:
\begin{lemma}
\label{comparison-form}
We have the comparison estimate on $\Omega$
\bea
-u_s\leq \varepsilon(-\psi_{s,k}+ \Lambda)^{2n\over 2n+1}
\eea
where $\varepsilon$ is defined by $\varepsilon^{1+2 n}=\frac{(2n+1)^{2n}}{2^{4n}n^{4n} \gamma^2 } A_{s,k} $, and $ \Lambda = \Big( \frac{10 C_0 C_\theta n \varepsilon}{2n+1} \Big)^{2n+1}$.
\end{lemma}

\noindent
{\it Proof.} We have to show that the test function
$$\Phi = -\varepsilon (-\psi_{s,k}+\Lambda)^{\frac{2n}{2n+1}} - u_s$$
is $\leq 0$ on $\Omega$. This is clear if $\max \Psi$ is achieved at $\overline{\Omega}\backslash \Omega_s$. So assume $\max_\Omega \Phi = \Phi(x_{\max})$ for some $x_{\max}\in \Omega_s$. By the maximum principle we have $\ddbar \Phi|_{x_{\max}}\le 0$ and $d\Phi|_{x_{\max}} = 0$. Note that at $x_{\max}\in \Omega_s\subset B(x_0, r_0)$, $d\Phi |_{x_{\max}}= 0$ implies that at $x_{\max}$ we have
 $$\varphi_i =  \frac{2n \varepsilon}{2n+1} (-\psi_{s,k} + \Lambda)^{-\frac{1}{2n+1}} (\psi_{s,k})_i -\epsilon' \bar z_i.$$
 Hence at $x_{\max}$ it follows that
\bea
|2 Re(G^{i\bar j} \varphi_i\bar\theta_j )| &= & \nonumber 2 |Re G^{i\bar j} \bar \theta_j (\frac{2n \varepsilon}{2n+1} (-\psi_{s,k} + \Lambda)^{-\frac{1}{2n+1}} (\psi_{s,k})_i -\epsilon' \bar z_i)| \\
&\le &2 C_0C_\theta \frac{2n\varepsilon}{2n+1} \Lambda^{-\frac{1}{2n+1}} \tr_G \omega_h + 2 \epsilon' r_0 C_\theta \tr_G \omega_h\label{eqn:6.5new}\\
&\le & \frac{2}{5} \tr_G \omega_h, \nonumber 
\eea
 if we choose the constant $\epsilon'>0$ small enough such that
 $$\omega_h \ge \frac{10\epsilon'}{n-1}(\tr_\omega \omega_h) \omega \mbox{ and } 2 \epsilon' r_0 C_\theta \le \frac 1 5. $$ In \neweqref{eqn:6.5new} we have applied \neweqref{eqn:Feb18} to bound $|D\psi_{s,k}|$ at $x_{\max}\in B(x_0,r_0)$. 
We define a linear operator $Lv = \Theta^{i\bar j} v_{i\bar j}$. The fact that $\Theta^{i\bar j}$ is positive definite implies that $L\Psi|_{x_{\max}}\le 0.$ We calculate at $x_{\max}$ as follows:
\bea
0& \ge& \nonumber L \Psi|_{x_{\max}} \\ 
& = & \nonumber \frac{2n\varepsilon}{2n+1} (-\psi_{s,k}+ \Lambda)^{-\frac{1}{2n+1}} \Theta^{i\bar j} (\psi_{s,k})_{i\bar j} + \frac{2\varepsilon n}{(2n+1)^2} (-\psi_{s,k}+ \Lambda)^{-\frac{2n+2}{2n+1}} \Theta^{i\bar j} (\psi_{s,k})_i (\psi_{s,k})_{\bar j} \\
 && \nonumber \quad  - 1 + \tr_{G} \omega_h -  \frac{\epsilon'}{n-1} ( \tr_{G} \omega \cdot \tr_{\omega} \omega_{{\mathbb C}^n} - \tr_{G} \omega_{{\mathbb C}^n} ) + 2 Re(G^{i\bar j} \varphi_i\bar\theta_j )\\
 & \ge & \nonumber \frac{2 n\varepsilon}{2n+1} (-\psi_{s,k}+ \Lambda)^{-\frac{1}{2n+1}} \Theta^{i\bar j} (\psi_{s,k})_{i\bar j}  - 1 + \tr_{G} [ \omega_h - \frac{2\epsilon'}{n-1} (\tr_\omega \omega_h) \omega ] - \frac{2}{5} \tr_G \omega_h\\
 & \ge & \nonumber \frac{4n^2\varepsilon}{2n+1} (-\psi_{s,k}+ \Lambda)^{-\frac{1}{2n+1}} (\det \Theta^{i\bar j})^{1/n} [\det ( \psi_{s,k})_{i\bar j}]^{1/n} - 1\\
 & \ge & \nonumber \frac{2 n^2\varepsilon \gamma^{1/n}}{2n+1} (-\psi_{s,k}+ \Lambda)^{-\frac{1}{2n+1}}\frac{(-u_s)^{1/2n}}{A_{s,k}^{1/2n}} - 1.
 \eea
 The term $2 Re(G^{i\bar j} \varphi_i\bar\theta_j )$ with first order derivatives of $\varphi$ is the only reason that we consider {\em real} Monge-Amp\`ere equation of $\psi_{s,k}$, instead of complex MA equations as in the sections above.  
 Here the first equality follows from the definition of $\Theta^{i\bar j}$ and the formula 
 $$- \Theta^{i\bar j} \varphi_{i\bar j} = - L\varphi = - 1 + \tr_{G} \omega_h + 2 Re( G^{i\bar j}\varphi_i \bar \theta_j) ,$$
 which is an easy consequence of the definition of $\tilde \omega$. In the fourth line, we applied the AG inequality. The last line follows from Lemma \ref{lemma 3.1} and \neweqref{eqn:Blocki}. By the choice of $\varepsilon$, this implies that $\Phi(x_{\max})\le 0$. Hence $\sup_\Omega \Phi\le 0$. Q.E.D.

\medskip

This means that we have now for the form type nonlinear equation the analogue of Lemma \ref{lm:comparison}. In fact, in the current case, the function $-\psi_{s,k} + \Lambda $ is uniformly bounded, hence the integral estimate \neweqref{eqn:alpha} still holds. We can then proceed in the same way as in Section \S 2 and apply Lemma \ref{general} to obtain the lower bound of $\varphi(x_0) = \inf_X \varphi$. We remark that the uniform $L^1(X,\omega^n)$ bound of $\varphi$ still holds in this setting. This is because the assumption that $\lambda[\tilde h_\varphi]\in \Gamma_1$ implies $\Delta_{\omega} \varphi + 2 Re(g^{i\bar j} \varphi_i \bar \theta_j)\ge - \tr_{\omega} \omega_h \ge -C$. A slight modification of  Lemma \ref{lemma 3} then yields the desired $L^1$ bound of $(-\varphi)$ given our normalization $\sup_X \varphi = 0$. The proof of Theorem \ref{thm:main4} is complete.

\section{The case of non-integrable almost complex structures}
\setcounter{equation}{0}
In this section, we generalize the comparison approach using {\em real} Monge-Amp\`ere equation as the auxiliary equation. 
Since $J$ is not integrable, there is no natural choice of local holomorphic coordinates on $X$. Instead, we will use local real coordinates and real Monge-Amp\`ere equation as the auxiliary equation. Let $r_0>0$ be a fixed number smaller than $1/10$ of the injectivity radius of $(X,\omega)$, viewed as a Riemannian manifold. Suppose $x_0\in X$ is a point where $\varphi$ achieves its minimum. We choose the normal coordinates $(x_1,\ldots, x_m)$ centered at $x_0$, where $m = 2n$. Without loss of generality, we may assume $B(x_0, 2r_0) = \{x| |x|< 2r_0\}$ is an Euclidean ball and in this ball
\bea
\label{real}
\frac 1 2 \delta_{ab}\le g_{ab}\le 2 \delta_{ab},\quad a, b = 1, 2,\ldots, m.
\eea 
Let $\eta>0$ be a small number to be determined. For any $s\in (0, s+0 = \eta r_0^2]$, we consider the function $u_s(x)$ defined on $B(x_0, 2r_0)$
\bea
\label{us}
u_s(x) = \varphi(x) - \varphi(x_0) + \eta |x|^2 - s.
\eea
It is clear that $u_s\ge -s$ on $B(x_0, 2r_0)$. Define the sublevel set $\Omega_s = \{x\in B(x_0, 2r_0)| u_s(x)<0\}$. Note that $x_0\in \Omega_s$, so $\Omega_s$ is a nonempty open set. By the choice of $s$, we have $u_s(x)> \eta r_0^2 -s >0$ for any $x\in B(x_0,2r_0)\backslash B(x_0, r_0)$, hence $\Omega_s\subset B(x_0,r_0)$. We solve the following real Monge-Amp\`ere equation on $B(x_0, 2r_0)$
\bea
\label{MongeAmpere}
\det \Big(\frac{\partial^2 \psi_{s,k}}{\partial x_a \partial x_b}\Big) = \frac{\tau_k(-u_s)}{A_{s,k}} e^{2nF},
\eea 
and $\psi_{s,k}$ is convex such that $\psi_{s,k} = 0$ on $\partial B(x_0, 2r_0)$. Here $A_{s,k}>0$ is chosen so that $\int_{B(x_0, 2r_0)} \det \Big( \frac{\partial^2 \psi_{s,k}}{\partial x_a \partial x_b}\Big) dx = 1$, i.e. 
$A_{s,k} = \int_{B(x_0, 2r_0)} \tau_k(-u_s) e^{2nF} dx \to \int_{\Omega_s} (-u_s) e^{2nF} dx=: A_s.$

\smallskip

We observe that the estimates \neweqref{eqn:Feb18} in Lemma \ref{lemma 11} hold for this $\psi_{s,k}$, in particular, $\sup_{B(x_0,r_0)} | D \psi_{s,k}|\le C_0$ for some $C_0 = C(n)$.

Define a test function $\Phi$ on $B(x_0, 2r_0)$ by 
\bea
\label{PHI}
\Phi(x) = - \varepsilon (-\psi_{s,k} (x) + \Lambda)^{\frac{2n}{2n+1}} - u_s(x),
\eea
where 
$$\varepsilon = \Big(\frac{4n+2}{n^2} \Big)^{\frac{2n}{2n+1}} {\gamma^{-2/(2n+1)}}  A_{s,k}^{1/(2n+1)}, \quad \Lambda \mbox{ to be determined}.$$
We claim that $\Phi\le 0$ in ${B(x_0,2r_0)}$. We may assume the maximum of $\Phi$ is obtained at $x_{\max}\in \Omega_s$, otherwise the claim follows easily. By maximum principle, at $x_{\max}$, $D\Phi = 0$ and $D^2 \Phi \le 0$. Define $L$ to be the linearized operator of $\varphi\mapsto \log f(\lambda[h_\varphi])$. Take a local unitary frame (w.r.t. $\omega$) $\{e_1,\ldots, e_n\}$ of $T^{(1,0)}X$ at $x_{\max}$ and let $\{\theta^1,\ldots \theta^n\}$ be the co-frame. Write $g_{i\bar j} = g(e_o,\bar e_j)$. Then 
$$\omega = g_{i\bar j} \sqrt{-1} \theta^i\wedge \bar \theta^j, \mbox{ and } g_{i\bar j} = \delta_{ij} \mbox{ at }x_{\max}.$$
By a unitary change of $e_i$, we may assume that $\omega_\varphi = \omega + \ddbar \varphi$ is diagonal at $x_{\max}$ and $\omega_\varphi = \tilde g_{i\bar j} \sqrt{-1} \theta^i \wedge \bar \theta^j$, $\tilde g_{i\bar j} = g_{i \bar j} + \varphi_{i\bar j}$, where for a function $\varphi$,
$$\varphi_{i\bar j} := e_i\bar e_j \varphi - [e_i, \bar e_j]^{(1,1)}\varphi.$$ That is, $\varphi_{i\bar j}$ is decomposed as the sum a second order derivative and a linear term in the first order derivative $D\varphi$, such that the latter is bounded by $C_1 |D\varphi|$, for a uniform constant $C_1 = C_1(\omega,n)>0$. We now determine the constant  $\Lambda = (\frac{20 n C_0 C_1}{2n+1})^{2n+1}\varepsilon^{2n+1} = C(n, \omega, \gamma) A_{s,k}$. With these notations, $L\Phi = G^{i\bar j} \Phi_{i\bar j}$ at $x_{\max}$. Then we calculate at $x_{\max}$

\bea
0 &\ge & L\Phi\nonumber \\
& = & \frac{2n \varepsilon}{2n + 1} (-\psi_{s,k} + \Lambda)^{-\frac{1}{2n+1}} G^{i\bar j} (\psi_{s,k})_{i\bar j} +  \frac{2n \varepsilon}{(2n + 1)^2} (-\psi_{s,k} + \Lambda)^{-\frac{2n+2}{2n+1}} |\partial \psi_{s,k}|_{G}^2\nonumber \\
&& - G^{i\bar j} \varphi_{i\bar j} - 4 \eta \tr_G g \nonumber\\
&\ge & \frac{2n \varepsilon}{2n + 1} (-\psi_{s,k} + \Lambda)^{-\frac{1}{2n+1}} G^{i\bar j} D^2_{i\bar j} \psi_{s,k} - \frac{2n C_1 \varepsilon}{2n + 1} (-\psi_{s,k} + \Lambda)^{-\frac{1}{2n+1}} \tr_G g |D\psi_{s,k}| \nonumber\\
&& \nonumber + \tr_G g - 1- 4\eta \tr_G g\\
&\ge &\frac{2n^2 \varepsilon}{2n + 1} (-\psi_{s,k} + \Lambda)^{-\frac{1}{2n+1}} \frac{\gamma^{1/n}}{e^{F}}(\det [D^2_{i\bar j} \psi_{s,k}])^{1/n} - 1\nonumber \\
&\ge &\frac{2n^2 \varepsilon}{2n + 1} (-\psi_{s,k} + \Lambda)^{-\frac{1}{2n+1}} \frac{\gamma^{1/n}}{e^{F}}\{ 4^{-n} (\det [D^2_{x_a x_b} \psi_{s,k}])^{1/2}\}^{1/n} - 1\nonumber \\
&\ge &\frac{ n^2 \varepsilon}{4n + 2} (-\psi_{s,k} + \Lambda)^{-\frac{1}{2n+1}} \frac{\gamma^{1/n}}{e^{F}} (\frac{(-u_s) e^{2nF}}{A_{s,k}})^{1/2n} - 1\nonumber.
\eea
We have choosen $\eta = 1/10$ and the third inequality follows from the choice of $\Lambda$ and $\eta$, and the arithmetic geometric inequality. The fourth inequality follows from the comparison of the determinants of complex and real Hessians. This proves the claim that $\Phi\le 0$. We can then proceed in the same way as in Section \S 2 and apply Lemma \ref{general} to obtain the lower bound of $\varphi(x_0) = \inf_X \varphi$. We remark that the uniform $L^1(X,\omega^n)$ bound of $\varphi$ still holds in this setting by Lemma \ref{lemma 15} below. The proof of Theorem \ref{thm:main5} is complete.

\begin{lemma}\label{lemma 15}
Suppose $\lambda[h_\varphi]\in \Gamma\subset \Gamma_1$. If we normalize $\sup_X \varphi = 0$, then
$$\int_X (-\varphi) \omega^n \le C(n,\omega).$$
\end{lemma}
\noindent{\em Proof.} First, $\lambda[h_\varphi] \in \Gamma_1$ implies that $\Delta_\omega^C \varphi \ge -n$ where $\Delta^C_\omega\varphi = n\frac{\omega^{n-1} \wedge \ddbar \varphi}{\omega^n}$ is the canonical Laplacian of $\varphi$. Second, it follows from standard calculations (e.g. \cite{CTW}) that $\Delta_\omega^R \varphi = 2 \Delta^C_\omega \varphi + \tau(d\varphi)$ for some ``torsion vector field'' $\tau$ of $(X,\omega)$, where $\Delta^R_\omega\varphi$ is the usual Laplacian of $\omega$ viewed as a Riemannian metric. Hence $\Delta^R_\omega \varphi - \tau(d\varphi) \ge -n $. We can then apply the standard elliptic estimates (Theorem 8.18 in \cite{GT}) over a local covering of $X$ as in the proof of Lemma \ref{lemma 3} to get the $L^1$ estimate of $-\varphi$.

\section{The case of tame symplectic forms}
\setcounter{equation}{0}

The proof of Theorem \ref{thm:main7} is analogous to that of Theorem \ref{thm:main5}, so we only point out the main differences. Recall $m=2n$ is the real dimension. We denote the Riemannian metric associated to $\Omega$ by $g$, i.e. 
$$g(U,V)={1\over 2}(\Omega(U,JV)+\Omega(V,JU)).$$
Around a point $x_0\in X$ in which $\varphi$ achieves the minimum, we take a local {\em real} coordinates system $(x_1, \ldots, x_m)$ under which the metric $g$ satisfies (\ref{real}). 
The almost-K\"ahler condition implies that the Calabi-Yau equation $\tilde\o^n=e^{nF}\Omega^n$ is equivalent to the real equation $\det (\tilde g_{ab}) = e^{2n F} \det (g_{ab})$, where as before we use the letters $a, b,\ldots$ to denote the indices $\{1, 2,\ldots, m\}$. A straightforward computation using the condition that $\tilde \omega$ is almost K\"ahler shows that 
$$\tilde g^{ab} \tilde\Gamma_{ab}^c = - \frac{1}{2} \tilde g^{c d} J^b{}_a \frac{\partial J^a{}_b}{\partial x_d} - \tilde g^{ab} J^c{}_d \frac{\partial J^d{}_b}{\partial x_a},$$ which implies that the Laplace-Beltrami operator $\Delta_{\tilde g} u = \tilde g^{ab} \frac{\partial^2 u}{\partial x_a \partial x_b} - \tilde g^{ab} \tilde \Gamma_{ab}^c \frac{\partial u}{\partial x_c}$ is independent of the first order derivatives of $\tilde g$. 

We define the function $u_s$ as in (\ref{us}) and solve the auxiliary real Monge-Amp\`ere equation (\ref{MongeAmpere}). For parameters $\varepsilon$ and $\Lambda$, we consider the test function $\Phi$ as in (\ref{PHI}). The observation that first order derivatives of the unknown metric $\tilde g$ is not involved in $\Delta_{\tilde g}$ enables us to apply the maximum principle the linear operator $L\Phi: = \Delta_{\tilde g} \Phi$, and with appropriate choice of the constants $\varepsilon, \Lambda$, we can derive that $\Phi\le 0$. The rest of the proof is similar to other sections, so we omit the details. 

\section{Further developments}
\setcounter{equation}{0}

We conclude by mentioning several recent developments of interest. The first is an extension of the present method for $L^\infty$ estimates to parabolic equations \cite{CC2}. The second is an extension of the theory to H\"older estimates \cite{Cheng}, and the third is an application of the method to new bounds for diameters and Green's functions in K\"ahler geometry \cite{GPS}, resulting in a positive answer for some long standing questions on the K\"ahler-Ricci flow \cite{GPSS}. In a different direction, building on the progress on Green's functions in \cite{GPS}, the existence of Hermitian-Einstein metrics on stable reflexive sheaves on normal K\"ahler spaces in the sense of Grauert has been very recently established in \cite{Cetal}.

\bigskip
\bigskip
\noindent
{\bf Acknowledgements} The authors would like to thank Zhiqin Lu for very helpful discussions on the Laplacian of Hermitian metrics.

\bigskip

\noindent Department of Mathematics \& Computer Science, Rutgers University, Newark, NJ 07102 USA

\noindent bguo@rutgers.edu,

\medskip

\noindent Department of Mathematics, Columbia University, New York, NY 10027 USA

\noindent phong@math.columbia.edu

\end{document}